\documentclass{amsart}

\usepackage{upref}
\usepackage{amssymb}
\usepackage[newitem,newenum,neverdecrease]{paralist}
\usepackage[final]{graphicx}
\usepackage[abbrev]{amsrefs}
\makeatletter
\renewcommand{\BibLabel}{%
    \hfil\hyper@anchorstart{cite.\CurrentBib}\relax\thebib\hyper@anchorend.%
}
\DefineJournal{trmis}{}{Tr. Mat. Inst. Steklova}{}
\DefineJournal{actam1}{0001-5962}{Acta Math.}{Acta Mathematica}
\DefineJournal{dukmj}{}{Duke Math. J.}{}
\DefineJournal{funap1}{}{Funktsional. Anal.\ i Prilozhen.}{}
\DefineJournal{izvan4}{}{Izv. Akad. Nauk SSSR Ser. Mat.}{}
\DefineJournal{matsb}{}{Mat. Sb.}{}
\DefineJournal{crmasp}{}{C. R. Math. Acad. Sci. Paris}{}
\DefineJournal{intjac}{}{Internat. J. Algebra Comput.}{}

\DefinePublisher{ams}{Amer. Math. Soc.}{American Mathematical
  Society}{Providence, RI}

\newcommand{\d@sh}[2]{\unskip#1\thinspace#2\thinspace\ignorespaces}
\newcommand{\dash}{\d@sh\nobreak\ndash}
\newcommand{\Dash}{\d@sh\nobreak\mdash}
\makeatother

\DeclareMathOperator*{\Wr}{\text{\Large $\wr$}}
\DeclareMathOperator*{\WR}{\text{\LARGE $\wr$}}

\numberwithin{equation}{section}

\DeclareMathOperator{\Aut}{Aut}
\DeclareMathOperator{\IMG}{IMG}

\DeclareMathOperator{\St}{St}
\DeclareMathOperator{\Sym}{Sym}
\DeclareMathOperator{\rist}{rist}
\DeclareMathOperator{\supp}{supp}
\DeclareMathOperator{\SP}{sp}
\DeclareMathOperator{\SPAN}{span}

\newcommand{\A}{\mathcal{A}}
\newcommand{\C}{\mathbb{C}}
\newcommand{\G}{\mathcal{G}}
\newcommand{\R}{\mathbb{R}}
\newcommand{\Z}{\mathbb{Z}}

\newcommand{\ii}{\mathrm{i}}
\renewcommand{\H}{\mathcal{H}}

\newtheorem{theorem}{Theorem}[section]
\newtheorem{prop}[theorem]{Proposition}
\newtheorem{cor}[theorem]{Corollary}
\newtheorem{lemma}[theorem]{Lemma}

\theoremstyle{definition}
\newtheorem{defin}[theorem]{Definition}

\begin{document}

\title[The Spectral Problem, Substitutions and Iterated Monodromy]{The Spectral Problem, Substitutions\\ and Iterated
Monodromy}

\author{Rostislav Grigorchuk}
\thanks{The first author is partially supported by NSF grants DMS-0308985 and
DMS-0456185.}

\author{Dmytro Savchuk}

\author{Zoran {\v S}uni\'c}
\thanks{The third author is partially supported by NSF grant
DMS-0600975}

\renewcommand{\shortauthors}{R.~Grigorchuk, D.~Savchuk, Z. {\v S}uni\'c}

\address{Department of Mathematics, Texas A\&M University,
College Station, TX 77843-3368, USA}

\email[R.Grigorchuk]{grigorch@math.tamu.edu}

\email[D.Savchuk]{savchuk@math.tamu.edu}

\email[Z. {\v S}uni\'c]{sunic@math.tamu.edu}


\begin{abstract}
We provide a self-similar measure for the self-similar group $\G$
acting faithfully on the binary rooted tree, defined as the iterated
monodromy group of the quadratic polynomial $z^2+\ii$. We also provide
an $L$-presentation for $\G$ and calculations related to the
spectrum of the Markov operator on the Schreier graph of the action
of $G$ on the orbit of a point on the boundary of the binary rooted
tree.
\end{abstract}

\maketitle
\section*{Introduction}

It was observed recently that the class of self-similar groups
naturally appears in mathematics. The most recent examples come from
combinatorics and are related to one of the most famous combinatorial
problems known as Hanoi Towers Game~\cite{grigorchuk-s:cr-hanoi}.
Slightly older examples, due first of all to
Nekrashevych~\cite{nekrash:self-similar} (see also
\cite{bartholdi_gn:fractal}), are related to holomorphic dynamics and
random walks.

Self-similar groups can be defined as groups generated by all the
states of a (not necessarily finite) Mealey
automaton~\cite{gns00:automata}. Of particular interest and importance
is the case when the automaton is finite, since in that case the
obtained group is finitely generated.

Self-similar groups have many interesting and important properties.
The class of self-similar groups contains exotic examples, such as
groups of Burnside type or groups of intermediate growth, as well as
familiar examples, such as free groups or free products of finite
groups, that are well known and are regular objects of study in
combinatorial group theory.

One of the most remarkable discoveries in the recent years is the
observation, due to Nekrashevych, that the so-called iterated
monodromy groups (IMG), which can be related to any self-covering map,
belong to the class of self-similar groups and that, in the most
natural situations, there is an explicit procedure representing them
by finite automata.

Even in the case of quadratic maps over $\C$ one gets a rich
theory with wonderful applications both to holomorphic dynamics and
to group theory~\cites{nekrash:self-similar,bartholdi-n:rabbit}.

Already the simplest examples of quadratic polynomials, such as
$z^2-1$ or $z^2+\ii$, show that the corresponding groups can be quite
complicated and can have extraordinary properties.

The group $\IMG(z^2-1)$ is called the Basilica group after the Julia
set of $z^2-1$ which (mildly) resembles the roof of the San Marco
Basilica in Venice (the top part of the Julia set is the roof and the
bottom part is its reflection in the water). Basilica group is torsion
free, of exponential growth, amenable but not elementary (and not even
subexponentially)
amenable~\cites{bartholdi_v:amenab,grigorch_z:basilica}, has trivial
Poisson boundary, is weakly branch, and has many other interesting
properties.

The main object of this article is the group $\IMG(z^2+\ii)$ (denoted by
$\G$ in the rest of the text), introduced
in~\cite{bartholdi_gn:fractal} and later studied by Bux and
Perez~\cite{bux_p:iter_monodromy}, who proved that $\G$ has
intermediate growth. This is not the first example of a self-similar
group of intermediate growth (the first examples were constructed
in~\cites{grigorch:burnside,grigorch:degrees}), but it is the first
example of a group of intermediate growth that naturally appears in
the area of applications of group theory.

We start with a quick introduction to the theory of self-similar
groups and, in particular, to iterated monodromy groups. We are
aiming for a self-contained treatment, which would make it possible
for the reader to understand the context of the paper completely
without reading other sources. In particular we give very detailed
calculation of the action of $\G= \IMG(z^2+\ii)$ on the binary
rooted tree. Then we show that the group $\G$ is a regular branch
group, thus presenting an example of a branch group which naturally
appears in holomorphic dynamics. The main body of the article is
devoted to the calculation of an $L$-presentation for $\G$ (i.e., a
presentation of a group by generators and relations which involves a
finite set of relators and their iterations by substitution).
Although it is known that $L$-presentations are quite common for
groups of branch type the number of examples in which explicit
computation is possible is rather small.

The presence of $L$-presentations is important from different points
of view. Such presentations are at the first level of complexity
after the finite presentations and quite often provide the simplest
way to describe a group that is not finitely presented ($\G$ is not
finitely presented~\cite{nekrashevych:cantor}). Further, such
presentations can be used to embed a group into a finitely presented
group in a way that preserves many properties of the original group.
We use the obtained $L$-presentation of $\G$ to embed $\G$ into a
finitely presented group with $4$ generators and 10 relators,
which is amenable but not elementary amenable (this approach has
been used for the first time in~\cite{grigorch:example}).

The rest of the article deals with finding a self-similar measure on
$\G$. The notion of a self-similar measure was introduced by
Kaimanovich in~\cite{kaiman:munchhausen}, who extends some ideas (in
particular the idea of self-similarity of a random walk) that appeared
in the work of Bartholdi and Vir\'ag~\cite{bartholdi_v:amenab}.

The self-similar measure is closely related to the problem of
computation of the spectrum of a Hecke type operator that can be
related to any group acting on a rooted tree and to the problem of
computation of the spectrum of the discrete Laplace operator (or,
what is almost the same, the Markov operator) on the boundary
Schreier graph of a group (i.e., the graph of the action of the
group on the orbit of a point of the boundary). A general approach
to the spectral problem (which extends the ideas outlined
in~\cites{bartholdi_g:spectrum,grigorch_z:asympt_spectrum}) based on
a renormalization principle and leading to questions on amenability,
multidimensional dynamics and multiparametric self-similarity of
operators is described in~\cite{grigorch_nv:spectra}. Unfortunately,
the spectral problem is not solved yet in our situation. What we are
able to construct is a rational map on $\R^3$ whose proper invariant
set (shaped as a ``strange attractor'') gives the spectrum of the
Markov operator after intersection by a corresponding line. Here we
have a situation analogous to the case of Basilica
group~\cite{grigorch_z:basilica_sp}. Further efforts in the
description of the shape of the attractor (and hence of the
spectrum) are needed.

The Schreier graph in this case, viewed through a macroscope, has a
form of a dendrite and this is a reflection of a general fact relating
the geometry of Schreier graphs and Julia sets proved by
Nekrashevych~\cite{nekrash:self-similar}.

In any case, our computations allow us to find a self-similar
nondegenerate measure on $\G$, which gives a self-similar random walk
on the group. The study of asymptotic properties of such random walks
is a promising direction and will be one of our subsequent subjects of
investigation.

\section{Iterated monodromy groups}

The theory of iterated monodromy groups was developed mostly by
Nekrashevych. A very detailed exposition can be found in his
monograph~\cite{nekrash:self-similar}. Here we give a definition and
some basic properties of these groups.

Consider a path connected and locally path connected topological
space $M$. Let $M_1$ be an open and path connected subset of $M$ and
$f\colon M_1\to M$ be a $d$-fold covering map. Fix a base point
$t\in M$ and let $\pi_1(M,t)$ be the corresponding fundamental
group. The set of iterated preimages of $t$ under $f$ has a natural
structure of a $d$-ary rooted tree $T$. Namely, each point $s$ from
this set has exactly $d$ preimages $s_1,\dots,s_d$ and these
preimages are declared to be adjacent to $s$ in $T$. The $n$th level
of the tree $T$ consists of the $d^n$ points in the set $f^{-n}(t)$.
Note that although the intersection of $f^{-n}(t)$ and $f^{-m}(t)$
may be nonempty for $m\neq n$, we formally consider the set of
vertices of $T$ to be a disjoint union of the sets $f^{-n}(t)$, $n
\geq 0$.

There is a natural action of $\pi_1(M,t)$ on the tree $T$. Let
$\gamma\in\pi_1(M,t)$ be a loop based at $t$. For any point $s$ of
$f^{-n}(t)$, there is a unique preimage $\gamma_{[s]}$ of $\gamma$
under $f^{n}$ which starts at $s$ and ends at a point $s'$, which
also belongs to $f^{-n}(t)$. We define an action of $\gamma$ on $T$
by setting $\gamma(s)=s'$. This action induces a permutation of
$f^{-n}(t)$ because the preimages of $\gamma^{-1}$ starting at the
points of $f^{-n}(t)$ are defined uniquely as well. The group of all
permutations of $f^{-n}(t)$ induced by all elements of $\pi_1(M,t)$
is called the \emph{$n$th monodromy group} of $f$. If $\gamma(s)=s'$
then $\gamma(f(s))=f(s')$ since $f(\gamma_{[s]})=\gamma_{[f(s)]}$,
so $\gamma$ acts on $T$ by a tree automorphism.

The action of $\pi_1(M,t)$ on $T$ is not necessary faithful. Let $N$
be the kernel of this action.

\begin{defin}
The group $\IMG(f)=\pi_1(M,t)/N$ is called the \emph{iterated
monodromy group} of $f$.
\end{defin}

It can be shown (see~\cite{nekrash:self-similar} for details) that, up
to isomorphism, $\IMG(f)$ does not depend on the choice of the base
point $t$.

In order to describe the automorphisms induced on $T$ by the loops
from $\pi_1(M,t)$ we need to come up with a ``coordinate system'' on
$T$. Let $X=\{0,1,\dots,d-1\}$ be a standard alphabet of cardinality
$d$. Then the set $X^*$ of all finite words over $X$ also has the
structure of a $d$-ary rooted tree, where $v$ is adjacent to $vx$,
for any $v\in X^*$ and $x\in X$.

The group $\Aut X^*$ of all automorphisms of $X^*$ has the structure
of an infinite iterated permutational wreath product $\wr_{i\geq
  1}\Sym(d)$ (because $\Aut X^*\cong \Aut X^*\wr_X\Sym(d)$,
where $\Sym(d)$ acts naturally on $X$ by permutations). This gives a
convenient way to express automorphisms from $\Aut X^*$ in the form
\begin{equation}\label{eqn_aut}
g=(g|_0,g|_1,\dots,g|_{d-1})\sigma_g,
\end{equation}
where $g|_0,g|_1,\dots,g|_{d-1}$ are automorphisms of the subtrees
of $X^*$ with roots at the vertices $0,1,\dots,d-1$ (these subtrees
are canonically identified with $X^*$) induced by $g$, and
$\sigma_g$ is the permutation of $X$ induced by $g$ (i.e.,
$\sigma_g(x)=g(x)$ \Dash the action of $g$ on $x\in X$). More
generally, for every $u\in X^*$ we define $g|_u$ to be the
automorphism of the subtree of $X^*$ rooted at $u$ (identified with
$X^*$) induced by $g$. The automorphism $g|_u$ is called the
\emph{section} of $g$ at $u$ and is uniquely determined by
$g(uw)=g(u)g|_u(w)$, for all $w\in X^*$.

The product of automorphisms written in form~\eqref{eqn_aut} is
performed in the following way. If
$h=(h|_0,h|_1,\dots,h|_{d-1})\sigma_h$ then
\[
gh=(g|_0h|_{\sigma_g(0)},\ldots,g|_{d-1}h|_{\sigma_g(d-1)})\sigma_g\sigma_h.
\]
By definition $gh(w)=h\bigl(g(w)\bigr)$.

\begin{defin}
  A group $G\leq\Aut X^*$ is called \emph{self-similar} if $g|_u\in G$
  for all $g\in G$ and $u\in X^*$.
\end{defin}

A convenient way to describe a particular finitely generated
self-similar group $G$ generated by automorphisms
$g_1,g_2,\dots,g_n$ is through a, so-called, \emph{wreath
recursion}. In this presentation we simply write the action of each
$g_i$ in the form
\[
g_i=\bigr(w_1(g_1,\dots,g_n),\dots,w_d(g_1,\dots,g_n)\bigl)\sigma_{g_i},
\]
where $w_i$, $i=1,\dots,n$, are words in the free group of rank~$n$.

Another language which describes self-similar groups is the language
of automaton groups (see the survey paper~\cite{gns00:automata} for
details).

\begin{defin}
  A Mealy \emph{automaton} is a tuple $(Q,X,\pi,\lambda)$, where $Q$
  is a set (a set of states), $X$ is a finite alphabet, $\pi\colon
  Q\times X\to Q$ is a transition function and $\lambda\colon Q\times
  X\to X$ is an output function. If the set of states $Q$ is finite
  the automaton is called \emph{finite}.
\end{defin}

One can think of an automaton as a sequential machine which, at each
moment of time, is in one of its states. Given a word $w\in X^*$ the
automaton acts on it as follows. It ``eats'' the first letter $x$ in
$w$ and depending on this letter and on the current state $q$ it
``spits out'' a new letter $\lambda(q,x)\in X$ and changes its state
to $\pi(q,x)$. The new state then handles the rest of word $w$ in
the same fashion. Thus the map $\lambda$ can be extended to
$\lambda\colon Q\times X^*\to X^*$ \Dash we just feed the automaton
with letters of $u\in X^*$ one by one. Each state $q$ of the
automaton defines a map, also denoted by $q$, from $X^*$ to itself
defined by $q(w)=\lambda(q,w)$. In the special case when, for all
$q\in Q$, the map $\lambda(q,\cdot)$ is a permutation of $X$ the map
$q\colon X^*\to X^*$ is invertible and hence, an automorphism of the
tree $X^*$. In this case the automaton is called \emph{invertible}.

\begin{defin}
  A group of automorphisms of $X^*$ generated by all the states of an
  invertible automaton $\A$ is called the \emph{automaton group}
  generated by $\A$.
\end{defin}

The class of automaton groups coincides with the class of self-similar
groups. Indeed, the action on $X^*$ of every element $g$ of a
self-similar group can be encoded by an automaton whose states are the
sections of $g$ on the words from $X^*$, transition and output
functions are derived from the representation~\eqref{eqn_aut}. Namely,
for each $u\in X^*$, set $\pi(g|_u,x)=g|_{ux}$ and
$\lambda(g|_u,x)=g|_u(x)$.

Important subclass of automaton groups consists of groups generated
by \emph{finite} automata. For example, we know that for groups in
this class the word problem is solvable.

\begin{figure}
\centering
\includegraphics[height=100pt]{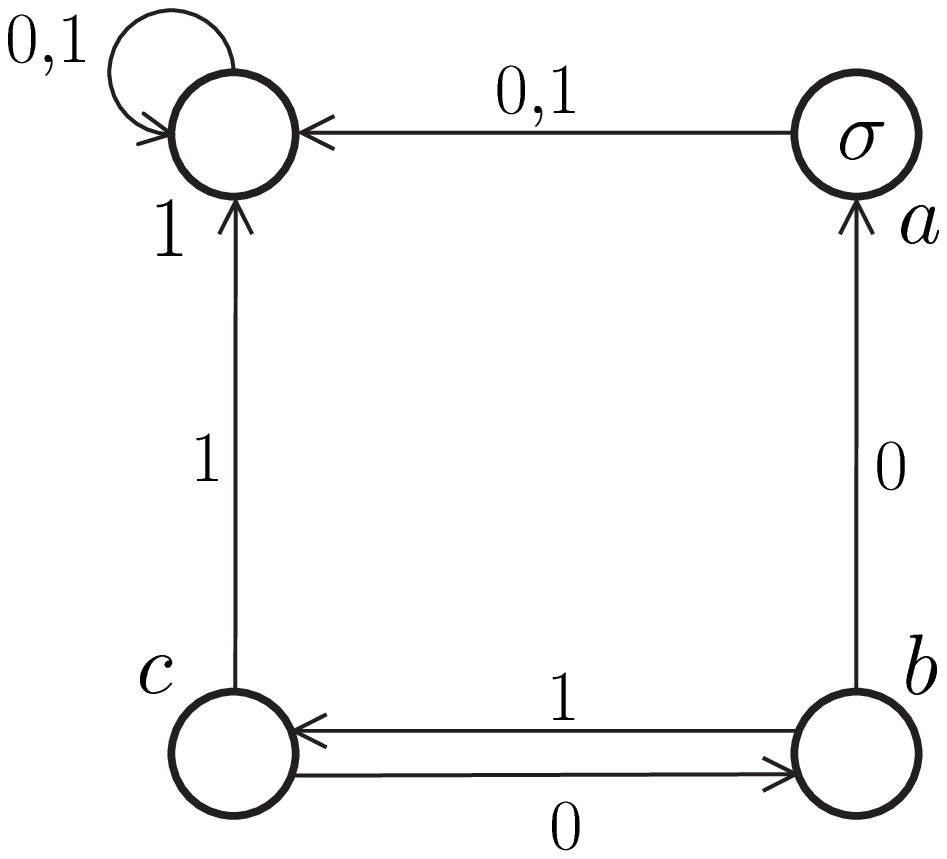}
\caption{Automaton generating group $\IMG(z^2+\ii)$\label{fig_aut}}
\end{figure}

A standard way to visualize automata is by so-called \emph{Moore
diagrams}. Such a diagram is an oriented graph where the set of
vertices is $Q$ and for every $q\in Q$, $x\in X$, there is an edge
from $q$ to $\pi(q,x)$ labeled by $\bigl(x,\lambda(q,x)\bigr)$. In
case of invertible automata it is common to label states by the
corresponding permutations of $X$ and leave only the first
coordinate on the edge labels. An example of a Moore diagram is
presented in Figure~\ref{fig_aut}.

We go back now to iterated monodromy groups and construct an
isomorphism $\Lambda\colon X^*\to T$ such that the induced action of
$\pi_1(M,t)$ on $X^*$ becomes particularly nice (self-similar).

We construct $\Lambda$ level by level. Set $\Lambda(\varnothing)=t$.
For each vertex $v$ in $X^n$ we will construct a path $l_v$ in $M$
joining $t$ to one of its preimages $s_v$ from $f^{-n}(t)$ and
define $\Lambda(v)=s_v$. Choose arbitrarily $d$ paths
$l_0,\dots,l_{d-1}$ in $M$ connecting $t$ to its $d$ preimages in
$f^{-1}(t)$ and, for $x\in X$, define $\Lambda(x)$ to be the end of
the path $l_x$. Now assume we have already defined $\Lambda(v)$ and
corresponding paths $l_v$ for all $v\in X^m$, $m\leq n$ and
$\Lambda$ is an isomorphism between the first $n$ levels of $X^*$
and $T$ such that, for all vertices $v$ on the first $n$ levels,
$\Lambda(v)$ is the endpoint of $\ell_v$. For any word $xu\in
X^{n+1}$ with $x\in X$ and $u\in X^n$ define
\[
l_{xu}= l_u f^{-n}_{[\Lambda(u)]}(l_x),
\]
where $f^{-n}_{[\Lambda(u)]}(l_x)$ is the unique preimage of the
path $l_x$ under $f^{n}$ starting at the vertex $\Lambda(u)$
(composition of paths is performed from left to right, i.e., the
path on the left is traversed first). Define $\Lambda(xu)$ to be the
end of the path $l_{xu}$.

\begin{figure}[b]
\centering
\includegraphics[height=100pt]{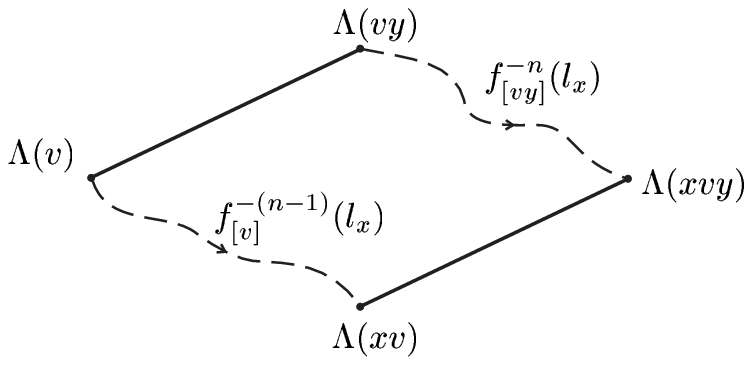}
\caption{Isomorphism $\Lambda$ between $T$ and
$X^*$\label{fig_img_isom}}
\end{figure}

In order to prove that $\Lambda$ is an isomorphism of trees we need
to show that $f\bigl(\Lambda(xvy)\bigr)=\Lambda(xv)$, for all $x,y\in X$ and
$v\in X^*$. Indeed,
\[
f(l_{xvy})=f(l_{vy})f(f^{-n}_{[\Lambda(vy)]}(l_x))=f(l_{vy})f^{-(n-1)}_{[\Lambda(v)]}(l_x).
\]
By definition, $f^{-(n-1)}_{[\Lambda(v)]}(l_x)$ is a path going from
$\Lambda(v)$ to $\Lambda(xv)$, so the end $\Lambda(xvy)$ of the path
$l_{xvy}$ is mapped to $\Lambda(xv)$ under $f$. Abusing the
notation, we often identify the trees $T$ and $X^*$ and write $v$
for $\Lambda(v)$ (see Figure~\ref{fig_img_isom}, where solid lines
represent edges in the tree $T$ and dashed lines represent paths in
$M$).

\begin{defin}
  The action of $\IMG(f)$ on $X^*$ induced by the isomorphism
  $\Lambda$ is called the \emph{standard action} of $\IMG(f)$.
\end{defin}

The tree isomorphism $\Lambda$ allows us to compute iterated
monodromy groups using the language of self-similar
groups~\cite{nekrash:self-similar}. We provide the details here in
order to keep the paper relatively self-contained and to help the
understanding of the computations that follow. Recall, that for any
loop $\gamma$ based at $t$ and any $u\in f^{-n}(t)$ we denote by
$\gamma_{[u]}$ the unique preimage of $\gamma$ under $f^{n}$ that
starts at the point $u$. Similarly, $f^{-n}_{[u]}(l_x)$ denotes the
unique preimage of the path $l_x$ starting at $u$.

\begin{theorem}
  The standard action of $\IMG(f)$ is self-similar. More precisely,
  the section $\gamma|_x$ of $\gamma\in \IMG(f)$ at $x\in X$ is given
  by
\begin{equation}\label{eqn_section}
\gamma|_x = l_x \gamma_{[x]}(l_{\gamma(x)})^{-1}.
\end{equation}
\end{theorem}

\begin{figure}
\centering
\includegraphics[height=140pt]{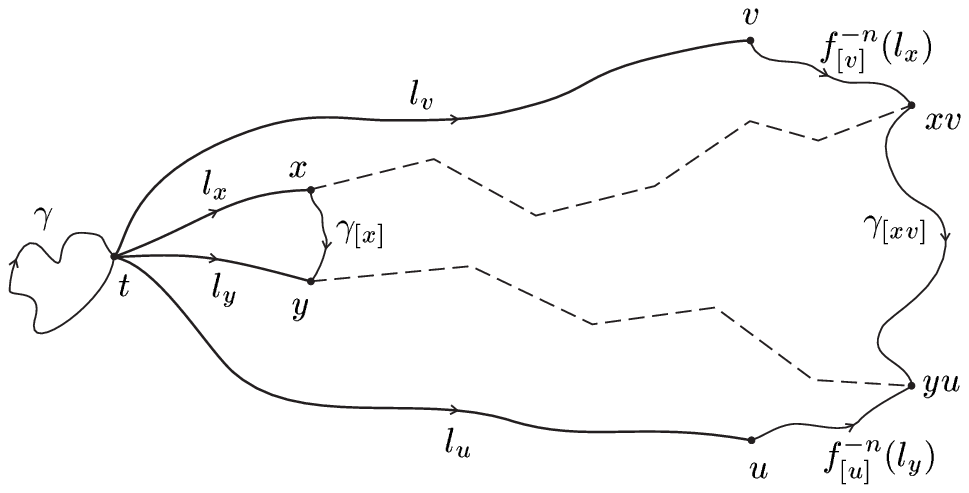}
\caption{Self-similar action of Iterated Monodromy
  Group\label{fig_img_selfsim}}
\end{figure}

\begin{proof}
  Let $v\in X^n$ be an arbitrary word and suppose $\gamma(xv)=yu$, for
  $y\in X$ and $u\in X^n$. Then vertices $v$ and $u$ are connected by
  the path
\[
p=f^{-n}_{[v]}(l_x) \cdot\gamma_{[xv]}\cdot
\bigl(f^{-n}_{[u]}(l_y)\bigr)^{-1},
\]
which goes through the vertices $v\to xv\to yu\to u$ (see
Figure~\ref{fig_img_selfsim}, where solid curves represent paths in
$M$ and dashed lines represent paths in the tree $X^*$)). We have
\[
f^n(p)= l_x \gamma_{[x]} l_y^{-1}.
\]
Thus the loop $\ell=l_x \gamma_{[x]} l_y^{-1}$ based at $t$
represents the element of $\IMG(f)$ which moves $v$ to $u$. The loop
$\ell$ is independent of $v$ (and $u$). Thus we have $\gamma|_x =
l_x \gamma_{[x]} l_y^{-1}$.
\end{proof}

We are now ready to compute the standard action of $\IMG(z^2+\ii)$
on $\{0,1\}^*$. The only critical point of this map is $z=0$, which
is preperiodic:
\[
0 \xrightarrow{f} \ii \xrightarrow{f} (\ii-1)
\overset{f}{\rightleftarrows} -\ii.
\]
and, hence, the postcritical set of $f$ is $\{\ii,\ii-1,-\ii\}$. Therefore
the restriction of $f$ on $M_1= \C \setminus\{\ii,\ii-1,-\ii,0\}$ is a
$2$-fold covering of $M= \C \setminus \{\ii,\ii-1,-\ii\}$.

\begin{figure}
\centering
\begin{tabular}{cc}
\includegraphics[height=140pt]{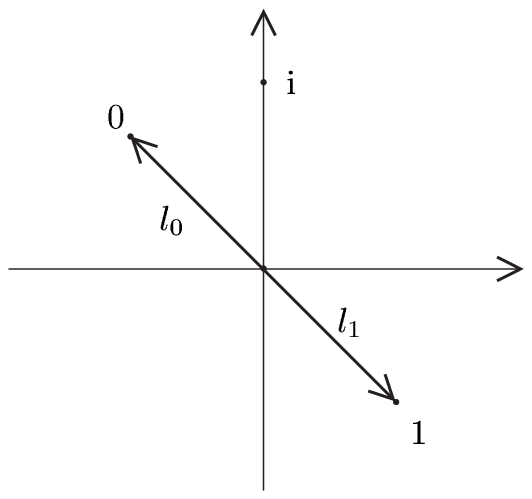}&
\includegraphics[height=140pt]{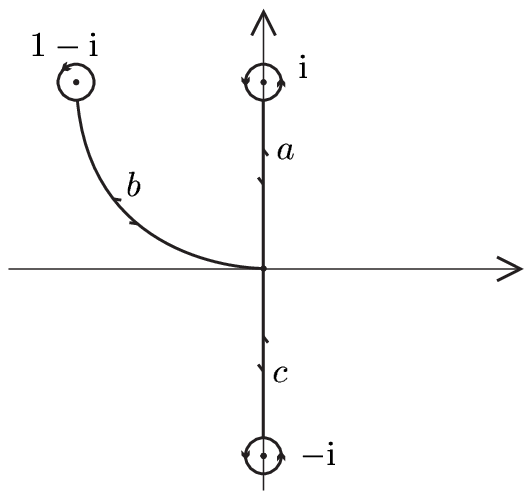}\\
(a)& (b)\\
\end{tabular}
\caption{(a)  Paths connecting $t=0$ to its preimages\hspace{2cm}
(b)  Generators of the fundamental group $\pi(M,t)$\label{fig_gens}}
\end{figure}

Set $t=0\in\C$ as the base point. It has two preimages $e^{\ii
3\pi/4}$ and $e^{\ii 7\pi/4}$ which are identified with the letters
$0$ and $1$, respectively (more precisely, we set $\Lambda(0)=e^{\ii
3\pi/4}$ and $\Lambda(1)=e^{\ii 7\pi/4}$). For the paths $l_0$ and
$l_1$ connecting $t$ to its preimages we choose the straight
segments shown in Figure~\ref{fig_gens}(a).

\begin{figure}
\centering
\includegraphics[width=\textwidth]{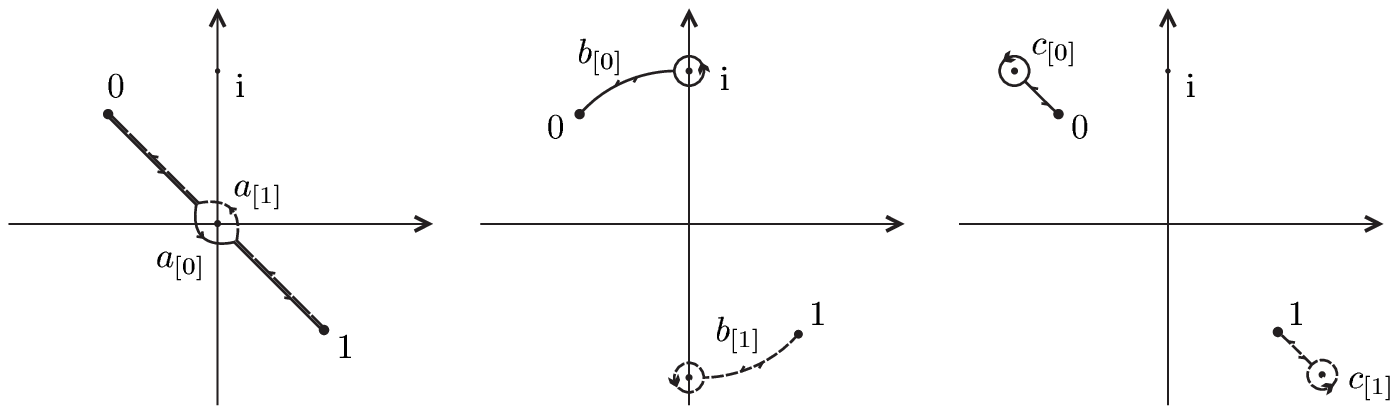}
\caption{Preimages of the generating loops\label{fig_gens_preim}}
\end{figure}

The fundamental group $\pi_1(M,t)$ is generated by the $3$ loops
$a,b,c$ shown in Figure~\ref{fig_gens}(b) going around $\ii$, $-\ii$
and $\ii-1$ respectively. Each of these loops has two preimages
$a_{[i]}$, $b_{[i]}$ and $c_{[i]}$, $i=0,1$, shown in
Figure~\ref{fig_gens_preim}.

According to formula~\eqref{eqn_section} and Figures~\ref{fig_gens}
and~\ref{fig_gens_preim} the sections of the generators $a$, $b$,
$c$ at $0$ and $1$ satisfy:
\[
\begin{aligned}
a|_0 &= l_0 a_{[0]} l_1^{-1} = 1,&\quad
a|_1 &= l_1 a_{[1]} l_0^{-1} = 1,\\
b|_0 &= l_0 b_{[0]} l_0^{-1} = a,&
b|_1 &= l_1 b_{[1]} l_1^{-1} = c,\\
c|_0 &= l_0 c_{[0]} l_0^{-1} = b,&
c|_1 &= l_1 c_{[1]} l_1^{-1} = 1,
\end{aligned}
\]
where $1$ denotes the trivial loop at $t$, which represents the
identity element of $\IMG(z^2+\ii)$.

Since $a$ permutes the elements of $f^{-1}(t)$, while $b$ and $c$ do
not, we obtain the following wreath recursion for the generators of
$\IMG(z^2+\ii)$
\begin{equation}\label{eqn_defin}
a=(1,1)\sigma,\quad
b=(a,c),\quad
c=(b,1),
\end{equation}
where $\sigma$ is the nontrivial transposition in $\Sym(2)$.

These relations show that the set of all sections of the generators
$a,b,c$ of $\G=\IMG(z^2+\ii)$ is $\{1,a,b,c\}$ and that the group $\G$
is generated by the states of the finite automaton shown in
Figure~\ref{fig_aut}.

We now say a few words about the relation between the dynamics of
the map $z\mapsto z^2+\ii$ and the combinatorial properties of the
action of $\G$ on the tree~$T$.

Recall that if a group $G$ acts on a set $Y$ then the \emph{Schreier
graph} of this action is an oriented graph, whose set of vertices is
$Y$ and there is an edge from $y\in Y$ to $z \in Y$ labeled by $g\in
G$ if and only if $g(y)=z$. It is convenient sometimes to forget
about the labels and/or the orientation of the edges.

Every group acting on a rooted tree acts on each level of the tree.
The Schreier graphs of such actions are of particular interest,
since in many situations (such as the one we are in) they can be
used to find the spectrum of the Markov operator on the boundary of
the tree.

\begin{figure}
\centering
\includegraphics[height=125pt]{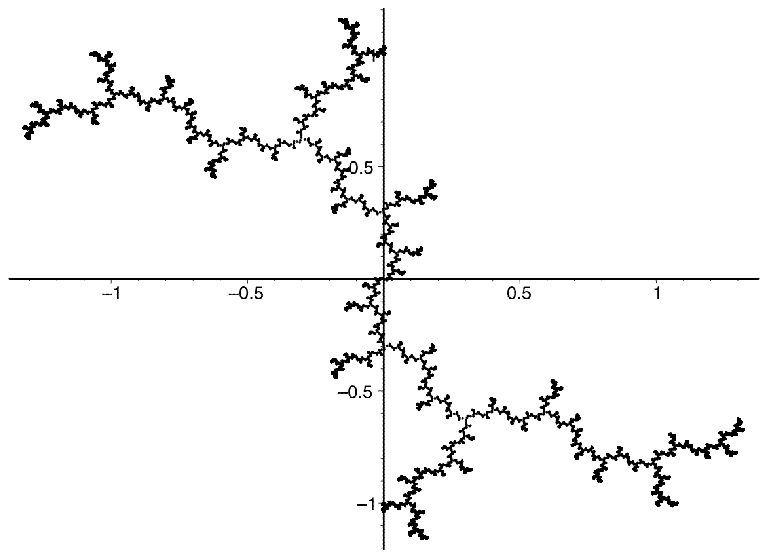}
\includegraphics[height=125pt]{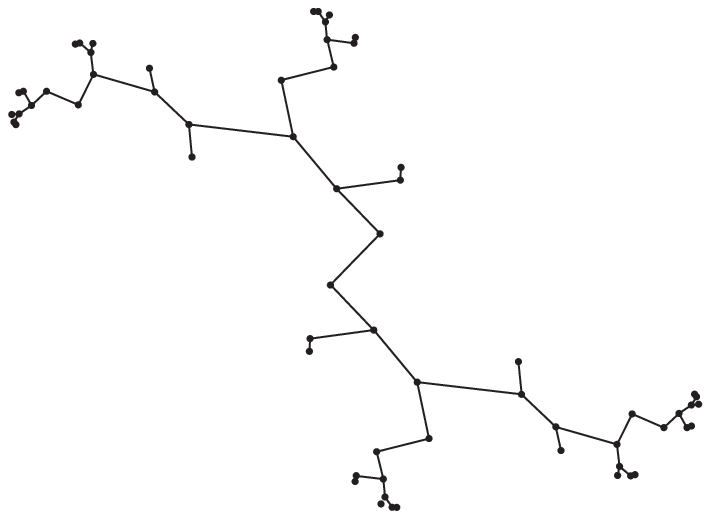}
\caption{Schreier graph of $\G$ and Julia set of
  $z^2+\ii$\label{fig_schreier}}
\end{figure}

Recent results of Nekrashevych~\cite{nekrash:self-similar} show that
the Schreier graphs of $\IMG(f)$ on the levels of the tree converge
to the Julia set of the map $f$. Therefore the structure of the
Julia set of $f$ provides understanding of the structure of the
Schreier graphs of $\IMG(f)$ (and vice versa). In our case the Julia
set of $z^2+\ii$ is the dendrite shown in the left half of
Figure~\ref{fig_schreier}. The right half of this figure displays
the Schreier graph of $\G$ on level $8$. The set of vertices of this
graph is just $f^{-8}(0)$ and the vertices are connected according
to the action of $\G$ (no loops are drawn though).

\section{Branch groups}

Another important class of subgroups of $\Aut X^*$ is the class of
branch groups~\cites{grigorch:jibranch,bar_gs:branch}. Here we give
basic definitions and prove that $IMG(z^2+\ii)$ is a regular branch
group. This shows that branch groups arise naturally in mathematics
(not just as a way to construct groups with unusual properties).

Let $G$ be a subgroup of $\Aut X^*$. Then for any vertex $v\in X^*$
one can define the subgroup of $G$ consisting of all the elements in
$G$ fixing all words in $X^*$ that do not have $v$ as a prefix. This
subgroup of $G$ is called the \emph{rigid stabilizer} of $v$ and it
is denoted by $\rist_G(v)$. Furthermore, the subgroup
\[
\rist_G(n)=\biggl\langle\bigcup_{v\in X^n}\rist_G(v)\biggr\rangle
\]
generated by the union of the rigid stabilizers of vertices at level
$n$, is called the \emph{rigid stabilizer of the $n$th level}.
Since elements of rigid stabilizers of different vertices on the
same level commute we have
\[
\rist_G(n)=\prod_{v\in X^n}\rist_G(v).
\]
Note that, if $G$ acts transitively on the levels of the tree, then
all rigid stabilizers of the vertices on a fixed level are conjugate
and, hence, isomorphic.

\begin{defin}
A group $G$ of tree automorphisms of $X^*$ that acts transitively on
the levels of $X^*$ is called a \emph{branch group} if all rigid
level stabilizers $\rist_G(n)$, $n \geq 0$, have finite index in
$G$. If all rigid stabilizers are nontrivial then $G$ is called a
\emph{weakly branch} group.
\end{defin}

It is often easier to prove that a given group belongs to a more
narrow class of regular (weakly) branch groups. Consider a
self-similar group $G$ and its normal subgroup $\St_G(1)$ consisting
of all elements in $G$ that stabilize the first level of $X^*$.
There is a natural embedding
\[
\Psi\colon \St_G(1)\hookrightarrow G\times G\times\dots\times G
\]
given by
\[
g\stackrel{\Psi}{\mapsto} (g|_0,g|_1,\dots,g|_{d-1}).
\]

\begin{defin}
Let $K,K_0,\dots,K_{d-1}$ be subgroups of a self-similar group $G$
acting on $X^*$. We say that $K$ \emph{geometrically contains} $K_0
\times \dots \times K_{d-1}$ and write
\[
K_0 \times \dots \times K_{d-1} \preceq K
\]
if $K_0 \times \dots \times K_{d-1} \leq \Psi(\St_G(1) \cap K)$.
\end{defin}

\begin{defin}
A group $G$ of tree automorphisms of $X^*$ that acts transitively on
the levels of the tree $X^*$ is called a \emph{regular weakly branch
group} over its normal subgroup $K$ if
\[
K \times \dots \times K \preceq K.
\]
If, in addition, the index of $K$ in $G$ is finite then $G$ is
called a \emph{regular branch group} over $K$.
\end{defin}

It can be shown that if $G$ is a regular (weakly) branch group than
it is a (weakly) branch group.

\begin{defin}
A self-similar group $G$ is called \emph{self-replicating} if, for
every vertex $u$ in $X^*$, the map $\varphi_u:G_u \to G$ given by
$\varphi_u(g)= g_u$ is onto (where $G_u$ is the stabilizer of the
vertex $u$ in $G$).
\end{defin}

Note that $\G=\IMG(z^2+\ii)$ is self-replicating. This is clear from
the equalities
\[
b=(a,c), \quad c = (b,1), \quad aba = (c,a), \quad aca = (1,b).
\]

Consider the normal subgroup $N$ of $\G$ defined by
\[
N = \langle[a,b],[b,c] \rangle^\G.
\]
By definition, $[g,h]=g^{-1}h^{-1}gh$ and $\langle\cdot\rangle^\G$
denotes normal closure in $\G$.

\begin{theorem}
  The group $\G$ is a regular branch group over $N$.
\end{theorem}

\begin{proof}
  First we observe that $N$ has finite index in $\G$. Direct
  computation shows that $a^2=b^2=c^2=(ac)^4=(ab)^8=(bc)^8=1$, so
  $\G/N$ is a homomorphic image of
\[
\langle a,b,c\mid a^2=b^2=c^2=(ac)^4=[a,b]=[b,c]=1\rangle \cong
\mathrm{C}_2\times \mathrm{D}_4,
\]
where $\mathrm{C}_2$ is the cyclic group of order 2 and $\mathrm{D}_4$
is the dihedral group of order 8.

Further, we have
\[
[b,c]=([a,b],1)\quad [c,b^a]=([b,c],1)
\]
Since $[c,b^a]=cb[b,a]cb[b,a]=[b,a]^{bc}[c,b][b,a]\in N$ we have that
$([a,b],1)$ and $([b,c],1)$ are elements in $N$. The fractalness of
$\G$ enables us to conjugate the sections in $([a,b],1)$ and
$([b,c],1)$ by arbitrary elements in $\G$ without leaving $N$. Thus we
get the inclusion $N\times 1 \preceq N$. The transitivity of $\G$ on
the first level then implies
\[
N \times N \preceq N.
\]

The level transitivity of $\G$ can be obtained almost for free. The
fact that $N$ is nontrivial along with the fact that $N\times N
\preceq N$ implies that $N$ is infinite, and hence so is $\G$. But one
can prove that a self-similar group acting on a binary rooted tree is
infinite if and only if it acts transitively on all levels.
\end{proof}

\section{\mathversion{bold}$L$-presentation}

The goal of this section is to prove the following result.

\begin{theorem}\label{thm_presentation}
The group $\G$ has the following $L$-presentation
\begin{multline}\label{eqn_pres}
\G\cong\bigl\langle a,b,c\bigm| \phi^n(a^2),\phi^n\bigl((ac)^4\bigr),
\phi^n([c,ab]^2),
\phi^n([c,bab]^2),\\ \phi^n([c,ababa]^2), \phi^n([c,ababab]^2),
\phi^n([c,bababab]^2), n\geq 0\bigr\rangle,
\end{multline}
where $\phi$ is the substitution defined on words in the free monoid
over the alphabet $\{a,b,c\}$ by
\[\phi\colon
\begin{cases}
a\to b,\\
b\to c,\\
c\to aba.\\
\end{cases}
\]
\end{theorem}

In order to prove Theorem~\ref{thm_presentation} we introduce some
notation and prove a few intermediate results.

The group
\[
\Gamma=\langle a,b,c\mid a^2,b^2,c^2,(ac)^4\rangle
\]
covers $\G$ (the relators of $\Gamma$ are relators of $\G$). The
action of $\G$ on the binary tree induces an action of the covering
group $\Gamma$ on the same tree, which is not faithful. Let $\Omega$
be the kernel of this action. Then, obviously, a set of generators
of $\Omega$ as a normal subgroup in $\Gamma$, together with the
relators in $\Gamma$ constitutes a presentation for $\G$.

The embedding $\G\hookrightarrow\G\mathbin{\wr}\Sym(2)$ induces a homomorphism
\[
\Psi\colon \Gamma\to\Gamma\mathbin{\wr}\Sym(2)
\]
defined by
\[
\Psi\colon\begin{cases}
a\mapsto(1,1)\sigma,\\
b\mapsto(a,c),\\
c\mapsto(b,1).
\end{cases}
\]
Indeed, the relators of $\Gamma$ are mapped to the trivial element
$(1,1)$ of $\Gamma\mathbin{\wr}\Sym(2)$:
\begin{align*}
\Psi(a^2)&=(1,1)\sigma(1,1)\sigma=(1,1),&
\Psi(b^2)&=(a,c)^2=(a^2,c^2)=(1,1),\\
\Psi(c^2)&=(b,1)^2=(b^2,1)=(1,1),&
\Psi((ac)^4)&=((1,b)\sigma)^4=(b^2,b^2)=(1,1).\\
\end{align*}

The homomorphism $\Psi$ induces homomorphisms $\Psi_n\colon
\Gamma\to\Gamma\mathbin{\wr}\bigl(\Wr_{i=1}^n\Sym(2)\bigr)$ (here
$\Wr_{i=1}^n\Sym(2)$ denotes the iterated permutational wreath
product) defined recursively by $\Psi_1=\Psi$ and
\[
\Psi_n\colon \Gamma\xrightarrow{\Psi_{n-1}}
\Gamma\mathbin{\wr}\biggl(\WR_{i=1}^{n-1}\Sym(2)\biggr)\xrightarrow{\Psi}
\bigl(\Gamma\mathbin{\wr}\Sym(2)\bigr)\mathbin{\wr} \biggl(\WR_{i=1}^{n-1}\Sym(2)\biggr)
=\Gamma\mathbin{\wr}\biggl(\WR_{i=1}^{n}\Sym(2)\biggr).
\]

If, for $g\in\Gamma$, we have $\Psi_n(g)=(g|_u,u\in X^n)\sigma_n$,
with $g|_u\in\Gamma$ and $\sigma_n\in \wr_{i=1}^n\Sym(2)$ we call
$g|_u$ the \emph{section} of $g$ at $u$.

For every $g\in\Gamma$ denote by $l(g)$ the length of the shortest
word in $a,b,c$ representing $g$ in $\Gamma$. The following lemma
shows that $\Gamma$ possesses the so called \emph{contraction}
property.

\begin{lemma}\label{lem_contr}
For every $g\in\Gamma$ and $u\in X^2$
\begin{equation}\label{eqn_contr}
l(g|_u)\leq \frac{l(g)+1}2.
\end{equation}
\end{lemma}

\begin{proof}
Observe that, because of the self-similarity, all generators satisfy
inequality~\eqref{eqn_contr}. Indeed,
\[
\Psi_2(a)=(1,1,1,1)(02)(13),\quad \Psi_2(b)=(1,1,b,1)(01),\quad
\Psi_2(c)=(a,c,1,1),
\]
where, for ease of notation, the vertices on the second level are
renamed by using the identifications $00 \leftrightarrow 0$, $01
\leftrightarrow 1$, $10 \leftrightarrow 2$, and $11 \leftrightarrow
3$.

All pairwise products of generators also satisfy
inequality~\eqref{eqn_contr}:
\begin{equation}
\begin{gathered}
\Psi_2(a^2)=1,\quad \Psi_2(b^2)=1,\quad \Psi_2(c^2)=1,\\
\begin{aligned}
\Psi_2(ab)&=\bigl(\Psi_1(c),\Psi_1(a)\bigr)\sigma=(b,1,1,1)(0213),\\
\Psi_2(ba)&=\bigl(\Psi_1(a),\Psi_1(c)\bigr)\sigma=(1,1,b,1)(0312),\\
\Psi_2(ac)&=\bigl(\Psi_1(1),\Psi_1(b)\bigr)\sigma=(1,1,a,c)(02)(13),\\
\Psi_2(ca)&=\bigl(\Psi_1(b),\Psi_1(1)\bigr)\sigma=(a,c,1,1)(02)(13),\\
\Psi_2(bc)&=\bigl(\Psi_1(ab),\Psi_1(c)\bigr)=(c,a,b,1)(01),\\
\Psi_2(cb)&=\bigl(\Psi_1(ba),\Psi_1(c)\bigr)=(a,c,b,1)(01).
\end{aligned}
\end{gathered}
\end{equation}
Any word $w$ in $a,b,c$ of length $n$ can be split into a product of
at most $(n+1)/2$ products of pairs of generators (if the length of
$w$ is odd one can pair the last letter in $w$ with $1$). Therefore
the sections of $w$ on the vertices of the second level are products
of at most $(n+1)/2$ letters. Thus the inequality~\eqref{eqn_contr}
holds for $w$ as well.
\end{proof}

Define an increasing sequence of subgroups of $\Gamma$ by
\[
\Omega_n=\ker\Psi_n.
\]

\begin{lemma}\label{lem_omega_union}
The kernel $\Omega$ of the canonical
epimorphism $\Gamma\to \G$ satisfies
\[
\Omega=\bigcup_{n\geq 1}\Omega_n
\]
\end{lemma}

\begin{proof}
Let $h$ be a word in $a,b,c$ of length at most $2^n+1$ representing
the trivial element in $\G$. Then, since for any words $u,v\in X^*$
\begin{equation}\label{eqn_restr_uv}
h|_{uv}=h|_u|_v,
\end{equation}
by Lemma~\ref{lem_contr} we obtain that all sections of $h$ have
length at most $1$ on the $2(n+1)$th level. Therefore they must be
trivial, because $h$ acts trivially on the tree. Hence,
$h\in\Omega_{2(n+1)}$.
\end{proof}

The last lemma reduces the problem of finding generators for $\Omega$
to finding generators for $\Omega_n$. We start from
$\Omega_1=\ker\Psi$ and, based on it, derive generators for
$\Omega_n$.

Let $H=\St_\Gamma(1)$ be the stabilizer of the first level of the
tree in $\Gamma$.

\begin{lemma}
The group $H$ has the following presentation
\[
H=\langle\beta,\delta,\gamma,\rho \mid
\beta^2=\delta^2=\gamma^2=\rho^2=(\rho\delta)^2=1\rangle,
\]
where $\beta=b$, $\delta=c$, $\gamma=aba$, $\rho=aca$.
\end{lemma}

\begin{proof}
  The index of $H$ in $\Gamma$ is $2$ and the coset representatives
  are $\{1, a\}$. The Reidemeister\dash Schreier procedure gives the
  above presentation.
\end{proof}

Obviously, each $\Omega_n$ is a subgroup of $H$. Therefore one can
restrict $\Psi$ to $H$. Since $H$ stabilizes the first level one can
think of $\Psi$ as a homomorphism $H\to\Gamma\times\Gamma$.  This map
is given by
\[\Psi\colon
\begin{cases}
\beta=b\to (a,c),\\
\gamma=aba\to (c,a),\\
\delta=c\to (b,1),\\
\rho=aca\to (1,b),\\
\end{cases}
\]
which mimics the corresponding embedding of the
generators $b$, $aba$, $c$, $aca$ of $\St_G(1)$ into $G\times G$.

Define the following words in $\Gamma$:
\[
\begin{aligned}
U_1&=(ba)^8,&\quad
U_2&=[c,ab]^2,&\quad
U_3&=[c,bab]^2,\\
U_4&=[c,ababa]^2,&
U_5&=[c,ababab]^2,&
U_6&=[c,bababab]^2.
\end{aligned}
\]

\begin{lemma}\label{lem_kernel}
$\Omega_1=\langle U_1, U_2, U_3, U_4, U_5, U_6\rangle^\Gamma$.
\end{lemma}

\begin{proof}
In order to find a generating set for $\Omega_1=\ker\Psi$ we first
describe $\Psi(H)$. Since $\Psi(\delta)=(b,1)$ and
$\Psi(\rho)=(1,b)$, we get
\[
B\times B \trianglelefteq \Psi(H),
\]
where $B=\langle b\rangle^\Gamma$. Furthermore, $\Psi(H)/(B\times
B)\cong\langle(a,c),(c,a)\rangle\cong \mathrm{D}_4$. Therefore
\[
\Psi(H)\cong(B\times B)\rtimes \mathrm{D}_4.
\]

Now we provide a presentation for $B$.

Define
\begin{alignat*}{4}
 \xi_1 &= b,   \qquad    & \xi_2 &= b^a, \qquad &
 \xi_3 &= b^c, \qquad    & \xi_4 &= b^{ca},\\
 \xi_5 &= b^{ac},        & \xi_6 &= b^{aca}, &
 \xi_7 &= b^{cac},       & \xi_8 &= b^{acac}.
\end{alignat*}
Since $\Gamma = C_2 * D_4$ where the cyclic group $C_2$ of order $2$
is generated by $b$ and the dihedral group $D_4$ of order $8$ is
generated by $a$ and $c$, it is clear that $B$ is generated by all
conjugates of $b$ by the elements in $D_4 = \langle a,c \rangle$.
Thus $\{\ \xi_i \mid i=1,\dots,8 \ \}$ is a generating set for $B$.
Moreover, it is clear that
\[
B=\langle\xi_i, i=1,\dots,8\mid \xi_i^2=1, i=1,\dots,8\rangle,
\]
i.e., $B$ is a free product of $8$ copies of the cyclic group of
order $2$ (indeed, none of the $b$'s in an expression of the form
$\xi_{i_1}\xi_{i_2}\dots\xi_{i_m}$ can be canceled in $\Gamma$ when
$i_j \neq i_{j+1}$ for $j=1,\dots,m-1$).

Therefore $B\times B$ is generated by $16$ elements, namely,
$\tilde\xi_i=(\xi_i,1)$ and $\hat\xi_i=(1,\xi_i)$ and is presented
by
\[
B\times B=\langle\tilde\xi_i, \hat\xi_j, i,j=1,\dots,8\mid
\tilde\xi_i^2, \hat\xi_j^2, [\tilde\xi_i,\hat\xi_j],
i,j=1,\dots,8\rangle.
\]

Now we compute the action of $\mathrm{D}_4$ generated by $x=(a,c)$ and
$y=(c,a)$ on $B\times B$.
\begin{equation}\label{eqn_rels}
\begin{aligned}
\tilde\xi_1^x&=(b,1)^{(a,c)}=(aba,1)=\tilde\xi_2,\\
\tilde\xi_2^x&=(aba,1)^{(a,c)}=(b,1)=\tilde\xi_1,\\
\tilde\xi_3^x&=(cbc,1)^{(a,c)}=(acbca,1)=\tilde\xi_4,\\
\tilde\xi_4^x&=(acbca,1)^{(a,c)}=(cbc,1)=\tilde\xi_3,\\
\tilde\xi_5^x&=(cabac,1)^{(a,c)}=(acabaca,1)=\tilde\xi_7,\\
\tilde\xi_6^x&=(cacbcac,1)^{(a,c)}=(acacbcaca,1)=\tilde\xi_8,\\
\tilde\xi_7^x&=(acabaca,1)^{(a,c)}=(cabac,1)=\tilde\xi_5,\\
\tilde\xi_8^x&=(acacbcaca,1)^{(a,c)}=(cacbcac,1)=\tilde\xi_6,\\
\tilde\xi_1^y&=(b,1)^{(c,a)}=(cbc,1)=\tilde\xi_3,\\
\tilde\xi_2^y&=(aba,1)^{(c,a)}=(cabac,1)=\tilde\xi_5,\\
\tilde\xi_3^y&=(cbc,1)^{(c,a)}=(b,1)=\tilde\xi_1,\\
\tilde\xi_4^y&=(acbca,1)^{(c,a)}=(cacbcac,1)=\tilde\xi_6,\\
\tilde\xi_5^y&=(cabac,1)^{(c,a)}=(aba,1)=\tilde\xi_2,\\
\tilde\xi_6^y&=(cacbcac,1)^{(c,a)}=(acbca,1)=\tilde\xi_4,\\
\tilde\xi_7^y&=(acabaca,1)^{(c,a)}=(cacabacac,1)=(acacbcaca,1)=\tilde\xi_8,\\
\tilde\xi_8^y&=(acacbcaca,1)^{(c,a)}=(cacacbcacac,1)=(acabaca,1)=\tilde\xi_7.
\end{aligned}
\end{equation}

The action on $\hat\xi_i$'s can be determined from the action on
$\tilde\xi_i$. Namely, if $\tilde\xi_i^x=\tilde\xi_p$ and
$\tilde\xi_i^y=\tilde\xi_q$, then
\begin{equation}
\label{eqn_rels_hat}
\hat\xi_i^x=\hat\xi_q\quad \text{and}\quad \hat\xi_i^y=\hat\xi_p.
\end{equation}

Now we can write down a presentation for $\Psi(H)$.
{\multlinegap0pt
\begin{multline*}
\Psi(H)=\langle\tilde\xi_i, \hat\xi_j, i,j=1,\ldots,8, x, y \mid
\tilde\xi_i^2=\hat\xi_j^2=[\tilde\xi_i,\hat\xi_j]=1,
i,j=1,\dots,8,\\
x^2=y^2=(xy)^4=1, \text{relations \eqref{eqn_rels}
and \eqref{eqn_rels_hat}}\rangle.
\end{multline*}}%

Note that relations~\eqref{eqn_rels} and~\eqref{eqn_rels_hat} show
that $\Psi(H)=\langle\tilde\xi_1,\hat\xi_1,x,y\rangle$. The kernel
of $\Psi$ is generated by the preimages of the relators of
$\Psi(H)$. We have
\[
\Psi(c)=\tilde\xi_1,\quad
\Psi(aca)=\hat\xi_1,\quad
\Psi(b)=x,\quad
\Psi(aba)=y.
\]
We can now start writing down the generators of $\ker\Psi$ as a
normal subgroup in $H$.
\[
x^2\to b^2=1\quad
y^2\to (aba)^2=1\quad
(xy)^4\to (baba)^4=(ba)^8=U_1
\]
We use the fact that $(ba)^8\in\ker\Psi$ to simplify the further
calculations.

Each of $\tilde\xi_i$'s and $\hat\xi_i$'s has the form
$\tilde\xi_1^{z_i}$ and $\hat\xi_1^{w_i}$, respectively, for some
$z_i,w_i\in\langle x,y\rangle$. The lifts of all elements from
$\langle x,y\rangle$ constitute the normal closure $\langle
b\rangle^{\mathrm{D}_8}$ of $b$ in
 $\mathrm{D}_8=\langle a,b\rangle$. Therefore the lifts of $\tilde\xi_i$'s will
have form
\[
\tilde\xi_i\to c^z,
\]
where $z$ runs over all elements of $\langle
b\rangle^{\mathrm{D}_8}$. In the same way the lifts of $\hat\xi_i$'s
look like
\[
\hat\xi_i\to(aca)^z=c^{az}=c^w,
\]
where $w$ runs over the complement of $\langle b\rangle^{\mathrm{D}_8}$ in
$\mathrm{D}_8$, which is just the coset $a\langle b\rangle^{\mathrm{D}_8}$.

Hence we get the remaining lifts of relators
\[
\begin{gathered}
\tilde\xi_i^2\to (c^z)^2=z^{-1}czz^{-1}cz=1,\\
\hat\xi_i^2\to (c^w)^2=w^{-1}cww^{-1}cw=1,\\
\bigl[\tilde\xi_i,\hat\xi_j\bigr],i,j=1,\dots,8 \to [c^z,c^w],
z\in\langle b\rangle^{\mathrm{D}_8}, w\in a\langle b\rangle^{\mathrm{D}_8}
\end{gathered}
\]

Now we would like to simplify the generators we obtained. First of
all, since $c^2=1$ we immediately get
\[
[c^z,c^w]=[c,c^{wz^{-1}}]^z
\]
so we can get rid of $z$ (because $wz^{-1}\in a\langle
b\rangle^{\mathrm{D}_8}$). Furthermore,
\[
[c,c^w]=cw^{-1}cw\cdot cw^{-1}cw=[c,w]^2
\]
We can discard $3$ more generators since
\[
[c,ba]^2=\bigl([c,ab]^{-2}\bigr)^c,\quad
[c,bababa]^2=\bigl([c,ababab]^{-2}\bigr)^c\quad
[c,a]^2=(ca)^4=1.
\]

Thus we get $5$ more generators for $\ker\Psi$:
\[
\begin{gathered}
{}[c,ab]^2=U_2,\quad
[c,bab]^2=U_3,\quad
[c,ababa]^2=U_4,\\
[c,ababab]^2=U_5,\quad
[c,bababab]^2=U_6.
\end{gathered}
\]

These generators, together with $U_1=(ba)^8$ generate $\ker\Psi$ as
a normal subgroup in $\Gamma$.
\end{proof}

\begin{lemma}\label{lem_omega_n}
$\Omega_{n}=\langle\phi^i(U_j), i=0,\dots,n-1,
j=1,\dots,6\rangle^\Gamma$
\end{lemma}

\begin{proof}
We will use induction on $n$. For $n=1$ the statement holds by
Lemma~\ref{lem_kernel}. Assume it is true for some fixed $n$.

By the definition of $\Omega_{n+1}$ we have
$\Psi(\Omega_{n+1})\leq\Omega_n\times\Omega_n$. We will show that
$\Psi(\Omega_{n+1})\geq\Omega_n\times\Omega_n$. Observe that
\begin{equation}\label{eqn_lifting}
\varphi_1\Bigl(\Psi\bigl(\phi^i(U_j)\bigr)\Bigr)=1
\end{equation}
in $\Gamma$ (recall that, for an element $h=(h|_0,h|_1)$ in $H$,
$\varphi_1(h)=h|_1$. Indeed, since
$\varphi_1\left(\Psi\bigl(\phi(\Gamma)\bigr)\right) \leq \langle
a,c\rangle=\mathrm{D}_4$ it's sufficient to check only that all
$U_i$'s are trivial in $\mathrm{D}_4$. But this is true since all
these words are squares of commutators and
$[\mathrm{D}_4,\mathrm{D}_4]\cong \Z/2\Z$.

Equation~\eqref{eqn_lifting} for $i=n+1$, together with the
inductive assumption yields $\Omega_n\times1 \leq
\Psi(\Omega_{n+1})$. Since $\Omega_{n+1}$ is normal in $\Gamma$
conjugation by $a$ yields $1\times\Omega_n \leq \Psi(\Omega_{n+1})$.
Therefore
\[
 \Psi(\Omega_{n+1})= \Omega_n\times\Omega_n.
\]
Equation~\eqref{eqn_lifting} also implies that
\[
\Psi(\phi^{n+1}(U_j))=(\phi^n(U_j),1),\quad
\Psi(\phi^{n+1}(U_j)^a)=(1,\phi^n(U_j)),
\]
i.e.,
\[
\Psi(\langle\phi^i(U_j), i=1,\dots,n,
j=1,\ldots,6\rangle^\Gamma)=\Omega_n\times\Omega_n.
\]
Therefore
\[
\begin{split}
\Omega_{n+1}&=\ker\Psi\cdot\langle\phi^i(U_j), i=1,\dots,n,
j=1,\dots,6\rangle^\Gamma\\
&=
\langle\phi^i(U_j), i=0,\dots,n, j=1,\dots,6\rangle^\Gamma.\qedhere
\end{split}
\]
\end{proof}

Lemmas~\ref{lem_omega_n} and~\ref{lem_omega_union} prove
Theorem~\ref{thm_presentation}. Since
$\phi\bigl((ac)^4\bigr)=(ba)^8=U_1$, $\phi(a^2)=b^2$, $\phi(b^2)=c^2$
and $\phi(c^2)=(aba)^2=1$ the presentation in~\eqref{eqn_pres} is
slightly simplified.

\begin{cor}
  The group $\G$ embeds into an amenable finitely presented group of
  exponential growth
\begin{multline*}
\tilde \G=\langle a,b,c,s\mid a^2, (ac)^4, [c,ab]^2,
[c,bab]^2, [c,ababa]^2, [c,ababab]^2, [c,bababab]^2,\\
a^s=b, b^s=c, c^s=aba\rangle,
\end{multline*}
which is an ascending HNN-extension of $\G$.
\end{cor}

\begin{proof}
By Theorem~\ref{thm_presentation} and the fact that
$\varphi_0(\Psi(\phi(u)))=u$ the substitution $\phi$ induces an
injective endomorphism of $\G$. Thus the HNN-extension construction
can be applied. Since $\tilde\G$ is an extension of the amenable
group $\G$ (see Section~\ref{sec_measures}) by the amenable group
$\Z$ generated by $s$, $\tilde\G$ is amenable. The growth of
$\tilde\G$ is exponential because it contains a free semigroup of
rank 2 (it follows from the HNN-extension construction that, for
example, $s$ and $sa$ generate such a semigroup).
\end{proof}

\section{Self-affine measures and amenability}
\label{sec_measures}

Although the amenability of $\G$ follows from the intermediate growth
of this group, which was established in~\cite{bux_p:iter_monodromy},
we present here a different approach based on the tools developed
in~\cites{bartholdi_v:amenab,kaiman:munchhausen}. More
precisely, we construct a particular self-affine measure on $\G$,
which proves the vanishing of the asymptotic entropy and, hence,
amenability.

Let $G$ be a self-similar group acting spherically transitively on a
$d$-ary tree. Consider a nondegenerate probability measure $\mu$ on
$G$ (the support of $\mu$ generates $G$). Then for any $x\in X$ one
can define a new probability measure $\mu|_x$ on $G$, which is
called the \emph{restriction} of $\mu$ on $x$. The details of the
definition and proofs of relevant statements are given
in~\cite{kaiman:munchhausen} and here we only give the basic idea.

We consider a right random walk $g_n=h_1h_2\ldots h_n$ on $G$
determined by $\mu$, i.e., $\{h_n\}$ is a sequence of independent
variables identically distributed according to the measure $\mu$. We
consider $\G$ as embedded in $\G\wr\Sym(2)$ and keep track of the
$x$th coordinate of the image of $g_n$ in $G\wr\Sym(d)$. Recall
that $h_n$ is an automorphism of $X^*$. For $x\in X$, $h_n(x)$
denotes the action of $h_n$ on $x$. Since
\[
g_{n+1}|_x=g_n|_x\cdot h_{n+1}|_{g_n(x)}
\]
the probability distribution of $g_{n+1}|_x$ is completely
determined by $g_n|_x$ and $g_n(x)$. Therefore the induced random
walk
\begin{equation}\label{eqn_rwidf}
(g_n|_x,g_n(x))
\end{equation}
on $G\times X$ is again a Markov chain. The last random walk is
called a \emph{random walk with internal degrees of freedom}. Since
$X$ is finite and $G$ acts transitively on $X$, the subset
$G\times\{x\} \subset G \times X$ is recurrent with respect
to~(\ref{eqn_rwidf}). Therefore one can consider the trace of the
random walk~\eqref{eqn_rwidf} on $G\times\{x\}$, which is also a
random walk. Finally, we define the measure $\mu|_x$ as the
transition law for the last random walk on $G\times\{x\}$ considered
as a copy of $G$.

There is a convenient way to compute $\mu|_x$ using the properties
of the random walk~\eqref{eqn_rwidf}. The random walk with internal
degrees of freedom is governed by the matrix
\[
M=(\mu_{xy})_{x,y\in X}
\]
whose entries $\mu_{xy}$ are subprobability measures on $G$ such that
$\mu_{xy}(h)$ is a transition probability of getting to the state
$(gh,y)$ from the state $(g,x)$.

With slight abuse of notation, we denote by $g$ the $\delta$-measure
concentrated at $g$. Then the matrix $M$ can be expressed as
\[
M = \sum_{g\in\supp\mu} \mu(g) M^g,
\]
where
\[
M^g_{xy}=
\begin{cases}
g|_x,&y=g(x),\\
0,&y\neq g(x),
\end{cases}
\]
The following theorem is proved in~\cite{kaiman:munchhausen}.

\begin{theorem}\label{thm_mu_x}
The measure $\mu|_x, x\in X$ can be expressed in
terms of the matrix $M$ as
\begin{equation}\label{eqn_mu_x}
\mu|_x=\mu_{xx}+M_{x\bar x}(1-M_{\bar{x}\bar{x}})^{-1}M_{\bar xx},
\end{equation}
where $M_{x\bar x}$ (resp., $M_{\bar xx})$ denotes the $x$th row
(column) of $M$ from which the entry $\mu_{xx}$ is removed, and
$M_{\bar{x}\bar{x}}$ is the matrix obtained from $M$ by removing the
$x$th row and the $x$th column.
\end{theorem}

One can define $\mu|_w$ for any $w=x_1x_2\dotsm x_n\in X^*$ by
\[
\mu_w=(\cdots(\mu|_{x_1})|_{x_2}\cdots)|_{x_n}.
\]

\begin{defin}
The nondegenerate probability measure $\mu$ on a self-similar group
$G$ is called \emph{self-affine} (\emph{self-similar}
in~\cite{kaiman:munchhausen}) if there is a word $w\in X^*$ such
that
\[
\mu|_w=\alpha e+(1-\alpha)\mu,
\]
where $0<\alpha<1$ and $e$ is the identity element in $G$.
\end{defin}

For simplicity, we write $\alpha$ instead of $\alpha e$.

\begin{theorem}[\cite{kaiman:munchhausen}]
If a self-similar group $G$ carries a self-affine nondegenerate
measure $\mu$ with finite entropy, then it is amenable.
\end{theorem}

In this section we construct such a measure on $\G$. Since this
measure should be nondegenerate and have finite entropy the most
natural place to look for it is the space $Q$ of positive convex
linear combinations of $\delta$-measures concentrated on the
generators $a$, $b$, and $c$, i.e., measure $\mu$ of the form
\[
\mu=xa+yb+zc,\quad x+y+z=1, \ x,y,z > 0.
\]
Suppose we want this measure to be self-affine with respect to $x\in
X$. By definition this means $\mu|_x=\alpha+(1-\alpha)\mu$ or,
equivalently,
\[
\mu=\frac{\mu|_x-\alpha}{1-\alpha}.
\]
Since $\mu(e)=0$ we get $\alpha=\mu|_x(e)$. Thus the measure $\mu$
is a fixed point of the transformation
\begin{equation}\label{eqn_Phi}
\Phi\colon \mu\mapsto\frac{\mu|_x-\mu|_x(e)}{1-\mu|_x(e)},
\end{equation}
which is defined in~\cite{kaiman:munchhausen} and used to prove
amenability of a family of groups generalizing Basilica group
$\IMG(z^2-1)$.

Let's compute $\mu|_0$ and the corresponding transformation $\Phi$
in the case of $\G$. The support of $\mu$ is $\{a,b,c\}$ and the
corresponding matrices $M^g$ are given by
\begin{equation}\label{eqn_oper_rec}
M^a=\begin{pmatrix}
0&e\\
e&0
\end{pmatrix},\quad
M^b=\begin{pmatrix}
a&0\\
0&c
\end{pmatrix},\quad
M^c=\begin{pmatrix}
b&0\\
0&e
\end{pmatrix}.
\end{equation}
Hence,
\[
M=xM^a+yM^b+zM^c=\begin{pmatrix}
ya+zb&x\\
x&yc+z
\end{pmatrix}.
\]
By Theorem~\ref{thm_mu_x}
\[
\mu|_0=(ya+zb)+x^2(1-yc-z)^{-1}.
\]
Since $c$ has order $2$ in $\G$ it is easy see that in the group
algebra $\R\G$
\[
(1-yc-z)^{-1}=\frac{y}{z^2-2z+1-y^2}\cdot
c-\frac{z-1}{z^2-2z+1-y^2}.
\]
Therefore
\[
\mu|_0=y\cdot a+z\cdot b+\frac{yx^2}{z^2-2z+1-y^2}\cdot
c-\frac{(z-1)x^2}{z^2-2z+1-y^2}
\]
and the transformation $\Phi$ takes the form
\[
\begin{split}
\Phi(xa+yb+zc)&=\frac{y\cdot a+z\cdot
b+ yx^2/(z^2-2z+1-y^2)\cdot
c}{1+(z-1)x^2/(z^2-2z+1-y^2)}\\
&= \frac{y(z^2-2z+1-y^2)}{z^2-2z+1-y^2+x^2z-x^2}\cdot a \\
&\qquad{}+ \frac{z(z^2-2z+1-y^2)}{z^2-2z+1-y^2+x^2z-x^2}\cdot b \\
&\qquad{}+ \frac{yx^2}{z^2-2z+1-y^2+x^2z-x^2}\cdot c.
\end{split}
\]

Now we are interested in a fixed point of the rational map
$F\colon \R^3\to \R^3$ induced by $\Phi$, which maps $(x,y,z)$ to the
coefficients of $\Phi(xa+yb+zc)$. Moreover, we are searching for
such a fixed point only in the invariant simplex $x+y+z=1$,
$x,y,z > 0$. Fortunately, there is such a fixed point. If
$\zeta\approx 0.4786202932$ is the unique real root of the polynomial
$Z^3-6Z^2+11Z-4$, then the point
\[
\bigl(\zeta,\zeta^2-4\zeta+2,-1+3\zeta-\zeta^2\bigr)
\]
is fixed under $F$, which produces a self-affine nondegenerate
probabilistic measure on $\G$ with finite support, proving
amenability of $\G$. This point is unique in the simplex of
nondegenerate measures. Indeed, from the equation
$F_1(x,y,1-x-y)=x$, where $F_1$ is the first coordinate of $F$, we
get
\begin{equation}\label{eqn_y}
y=\tfrac14(-x^2+x\pm x\sqrt{x^2-10x+9}).
\end{equation}
Substitution in $F_2(x,y,1-x-y)=y$ yields
\begin{equation}\label{eqn_root}
(-x^4+12x^3-39x^2+40x-12)\pm\sqrt{x^2-10x+9}(x^3-7x^2+12x-4)=0.
\end{equation}

\begin{figure}
\centering
\begin{tabular}{cc}
\includegraphics[height=120pt]{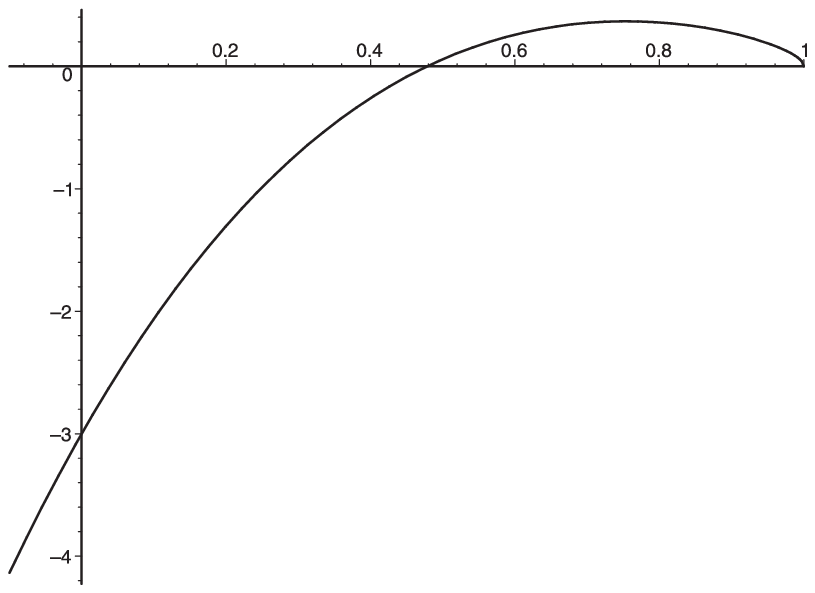}&
\includegraphics[height=120pt]{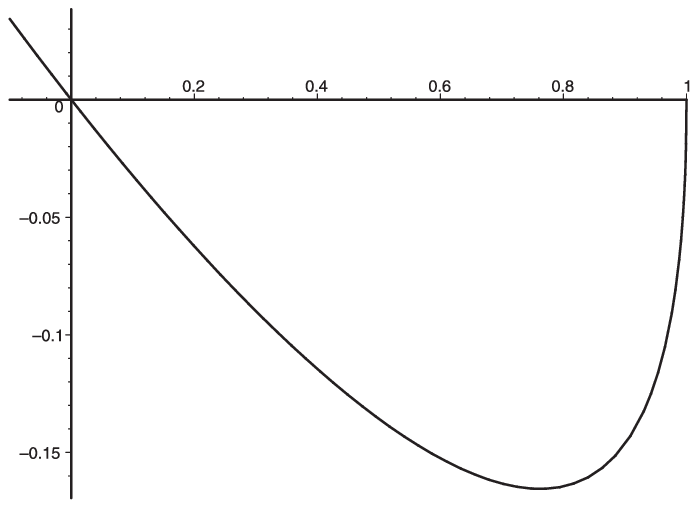}\\
(a)& (b)\\
\end{tabular}
\caption{Uniqueness of a self-affine measure\label{fig_root}}
\end{figure}

Moving the second summand to the righthand side and squaring both
sides produces the equation $x(x-1)(x^3-6x^2+11x-4)=0$, whose unique
real root on the interval $(0,1)$ is $\zeta$. The graphs of the two
functions in~\eqref{eqn_root} are shown in Figure~\ref{fig_root}:
(a) for ``plus'' and (b) for ``minus.'' The solution comes from
(a), so in~\eqref{eqn_y} ``plus'' should be used. It is a routine
to check that indeed $y=\zeta^2-4\zeta+2$.

\section{Spectral properties and Schur complement}
\label{sec_spectra}

Let $H$ be a Hilbert space and $M$ be an operator on $H$. Let
$H=H_0\oplus H_1$ and there are operators $A\in B(H_0)$, $D\in
B(H_1)$, $B\colon H_1\to H_0$ and $C\colon H_0\to H_1$, such that
the matrix of $M$ in the basis consisting of the bases of $H_0$ and
$H_1$ takes the form:
\[
M=
\begin{pmatrix}
A&B\\
C&D
\end{pmatrix}.
\]

The following fact is of folklore type.

\begin{prop}\label{prop_oper_inv}
  Let $D$ be invertible. The operator $M$ is invertible if and only if
  $S_1(M)=A-BD^{-1}C$ is invertible.
\end{prop}

The matrix $S_1(M)$ is called \emph{the first Schur complement} of
$M$.

\begin{proof}
Indeed, the matrix
\[
L=
\begin{pmatrix}
I_0&0\\
-D^{-1}C&D^{-1}
\end{pmatrix}
\]
is invertible. Therefore $M$ is invertible if and only if
\[
ML=
\begin{pmatrix}
A-BD^{-1}C&BD^{-1}\\
0&I_1
\end{pmatrix}
\]
is invertible, which is equivalent to the nonsingularity of
$S_1(M)$.
\end{proof}

In our case the action of $\G$ on the boundary $X^\omega$ (the set
of infinite sequences over $X$) of the tree $X^*$ induces a unitary
representation $\pi_g(f)(x)=f(g^{-1}x)$ of $\G$ in
$\H=B\bigl(L_2(X^\omega)\bigr)$. The Markov operator
$M=\tfrac13(\pi_a+\pi_b+\pi_c)$ corresponding to this unitary
representation plays an important role (we do not include inverse
elements because all generators are of order $2$). The usual method
to find the spectrum of $M$ for a self-similar group $G$ is to
approximate $M$ with finite dimensional operators arising from the
action of $\G$ on the levels of the tree $X^*$. For more on this
see~\cite{bartholdi_g:spectrum}. For simplicity we write $g$ for
$\pi_g$.

Let $\H_n$ be the subspace of $\H$ spanned by the $2^n$
characteristic functions $f_v$, $v\in X^n$, of the cylindrical sets,
corresponding to the $2^n$ vertices of the $n$th level. Then $\H_n$
is invariant under the action of $\G$ and $\H_n\subset \H_{n+1}$.
Also $\H_n$ can be naturally identified with $L_2(X^n)$. By
$\pi^{(n)}_g$ (or, with a slight abuse of notation, by $g_n$) we
denote the restriction of $\pi_g$ on $\H_n$. Then, for $n \geq 0$,
\[
M_n=\tfrac13(a_n+b_n+c_n)
\]
are finite-dimensional operators whose spectra converge to the
spectrum of $M$ in the sense
\[\
\SP(M)=\overline{\bigcup_{n\geq0}\SP(M_n)}.
\]

Moreover, if $P$ is the stabilizer of an infinite word from
$X^\omega$, then one can consider the Markov operator $M_{G/P}$ on
the Schreier graph of $G$ with respect to $P$. The following fact is
observed in~\cite{bartholdi_g:spectrum} and can be applied in the
case of $\G$.

\begin{theorem}
If $G$ is amenable then
\[
\SP(M_{G/P})=\SP(M).
\]
\end{theorem}

\begin{figure}[t]
\centering
\includegraphics*[height=200pt]{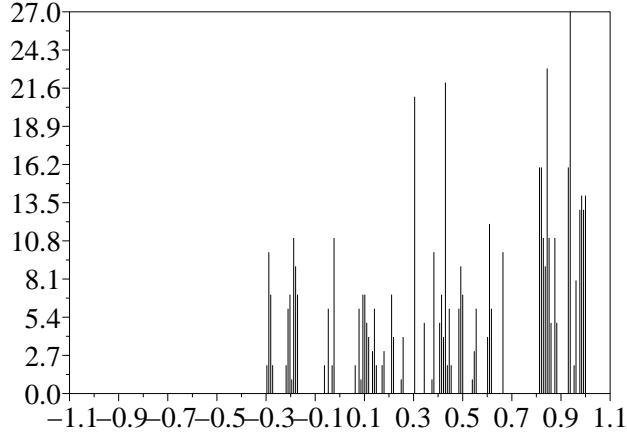}
\caption{Histogram of the spectrum on the 9th level\label{fig_spec}}
\end{figure}

Common practice for finding the spectrum of $M$, initiated
in~\cite{bartholdi_g:spectrum}, is to consider a pencil of operators
\[
\widetilde M(y,z,\lambda)=a+yb+zc-\lambda
\]
and find the set $\SP(y,z,\lambda)$ of points $(y,z,\lambda)$ such
that $\widetilde{M}(y,z,\lambda)$ is not invertible. Then the spectrum of
$M$ is just the intersection of $\SP(y,z,\lambda)$ with the line
$y=z=1$, shrunk by a factor of $\frac13$. We take $1$ as the
coefficient at $a$ to simplify the computation. Otherwise one can
divide it out (we restrict our attention to the case when $x,y,z$
are nonzero).

Let us consider the corresponding pencil $\tilde
M_n(y,z,\lambda)=a_n+yb_n+zc_n-\lambda$ and find its matrix in the
basis $\{f_v\!:\! v\in X^n\}$. The orthogonal subspaces
$H_n^{(i)}\!=\!\SPAN(f_{iv},{v\in X^{n-1}})$, $i=0,1$ span $H_n$ and
are naturally isomorphic to $H_{n-1}$. The self-similar structure of
$\G$ gives the following operator recursion (which coincides with
the recursion~\eqref{eqn_oper_rec})
\begin{equation}
a_n=\begin{pmatrix}
0&I_{n-1}\\
I_{n-1}&0
\end{pmatrix},\quad
b_{n}=\begin{pmatrix}
a_{n-1}&0\\
0&c_{n-1}
\end{pmatrix},\quad c_n=\begin{pmatrix}
b_{n-1}&0\\
0&I_{n-1}
\end{pmatrix},
\end{equation}
for $n>0$, where $I_{n-1}$ denotes the identity matrix of size
$2^{n-1}$. The matrices $a_0$, $b_0$ and $c_0$ are equal to the $1
\times 1$ matrix $[1]$. For any constant $r$, we write $r$ instead
of $rI_n$. Thus we have,
\[
\widetilde{M}_n(y,z,\lambda)=a_n+yb_n+zc_n-\lambda=\begin{pmatrix}
ya_{n-1}+zb_{n-1}-\lambda&1\\
1&yc_{n-1}+z-\lambda
\end{pmatrix}.
\]

By Proposition~\ref{prop_oper_inv} in case $yc_{n-1}+z-\lambda$ is
invertible the operator $\widetilde{M}_n(y,z,\lambda)$ is invertible if
and only if $S_1(\widetilde{M}_n(y,z,\lambda))$ is invertible. The
inverse of $yc_{n-1}+z-\lambda$ in $\R\G$ is
\[
\frac{y}{-y^2+z^2-2z\lambda+\lambda^2}\cdot c_{n-1}
+\frac{-z+\lambda}{-y^2+z^2-2z\lambda+\lambda^2}.
\]
Hence, $yc_{n-1}+z-\lambda$ is not invertible if and only if
$-y^2+z^2-2z\lambda+\lambda^2=({z-\lambda-y})\*({z-\lambda+y})=0$. Denote
the union of these 2 planes by $Z_1$. Note, that
$\widetilde{M}_n(y,z,\lambda)$ is not necessary singular at each point
of $Z_1$.

The first Schur complement of $\widetilde{M}$ is
{\multlinegap0pt
\begin{multline*}
S_1(\widetilde{M}_n(y,z,\lambda))\\
\begin{aligned}
&=
ya_{n-1}+zb_{n-1}-\lambda - (yc_{n-1}+z-\lambda)^{-1}\\
&=y\cdot a_{n-1}+z\cdot b_{n-1}
+\frac{y}{-y^2+z^2-2z\lambda+\lambda^2}\cdot c_{n-1}
+\frac{-z+\lambda}{-y^2+z^2-2z\lambda+\lambda^2}-\lambda
\end{aligned}.
\end{multline*}}%
If $y=0$ we get $S_1(\widetilde{M}_n(0,z,\lambda))=zb_{n-1}+\frac1{\lambda-z}-\lambda$ is not
invertible if and only if
\begin{multline*}
\det\begin{pmatrix}
1/(\lambda-z)-\lambda& z\\
z & 1/(\lambda-z)-\lambda
\end{pmatrix}\\
=
\frac1{(\lambda-z)^2}(1-(\lambda-z)(\lambda+z))(1-\lambda+z)(1+\lambda-z)=0.
\end{multline*}
Denote corresponding union of a hyperbola and two lines in $\R^3$ by
$Z_2$. Note that $Z_2\cap Z_1=\varnothing$.

If $y\neq0$ then
\[
\begin{split}
\frac1y S_1(\widetilde{M}_n(y,z,\lambda))
&=a_{n-1}+\frac zy\cdot b_{n-1}
+\frac{1}{-y^2+z^2-2z\lambda+\lambda^2}\cdot c_{n-1}\\
&\hspace{1.75in}{}-\frac{-\lambda y^2+\lambda
z^2-2z\lambda^2+\lambda^3+z-\lambda}{y(-y^2+z^2-2z\lambda+\lambda^2)}\\
&=\widetilde{M}_{n-1}(F(y,z,\lambda)),
\end{split}
\]
where $F\colon\R^3\to \R^3$ is the rational map defined by
\[
F\colon (y,z,\lambda)\to \biggl(\frac
zy,\frac{1}{-y^2+z^2-2z\lambda+\lambda^2}, \frac{-\lambda y^2+\lambda
z^2-2z\lambda^2+\lambda^3+z-\lambda}{y(-y^2+z^2-2z\lambda+\lambda^2)}\biggr).
\]

Therefore the set $\SP_n(y,z,\lambda)$ of points $(y,z,\lambda)$
where $\widetilde{M}_n(y,z,\lambda)$ is not invertible in this case
($y\neq 0)$ is a preimage under $F$ of the corresponding set
$\SP_{n-1}(y,z,\lambda)$. To summarize,
\begin{equation}
\label{eqn_sp_Z1}
Z_2\cup F^{-1}(\SP_{n-1}(y,z,\lambda))\subset
\SP_n(y,z,\lambda)\subset Z_1\cup Z_2\cup
F^{-1}(\SP_{n-1}(y,z,\lambda)).
\end{equation}

Since $\widetilde{M}_0(y,z,\lambda)=(1+y+z-\lambda)$ we have
\begin{equation}
\label{eqn_sp0}
\SP_0(y,z,\lambda)=\{(y,z,\lambda): 1+y+z-\lambda=0\}.
\end{equation}
Denote this plane by $P$.

Equations~\eqref{eqn_sp_Z1} and~\eqref{eqn_sp0} show that
\begin{equation}
\label{eqn_sp_incl}
F^{-n}(P)\cup\bigcup_{i=0}^{n-1}F^{-i}(Z_2)\subset
\SP_n(y,z,\lambda)\subset
F^{-n}(P)\cup\bigcup_{i=0}^{n-1}F^{-i}(Z_1\cup Z_2).
\end{equation}

Since $Z_2$ consists of points with $y=0$ every point
$(y,z,\lambda)$ from $F^{-1}(Z_2)$ must satisfy $z=0$. But the
preimages of all such points are empty. Hence,
$F^{-2}(Z_2)=\varnothing$ and $\bigcup_{i=0}^{n-1}F^{-i}(Z_2)=Z_2\cup
F^{-1}(Z_2)$. Denote the last subset by $Z_3$.

One can easily check that $P\subset F^{-1}(P)$ and, hence,
$F^{-n}(P)=\bigcup_{i=0}^n F^{-i}(P)$. Therefore
equation~(\ref{eqn_sp_incl}) transforms to
\[
Z_3\cup\bigcup_{i=0}^{n}F^{-i}(P)\subset \SP_n(y,z,\lambda)\subset
F^{-n}(P)\cup Z_3\cup\bigcup_{i=0}^{n-1}F^{-i}(P\cup Z_1).
\]
Thus, the spectrum of the operator $\widetilde{M}$ on $L_2(X^\omega)$
satisfies
\[
\overline{Z_3\cup\bigcup_{i=0}^\infty F^{-i}(P)}\subset \SP(\tilde
M(y,z,\lambda))=\overline{\bigcup_{i=0}^\infty
\SP_n(y,z,\lambda)}\subset\overline{Z_3\cup\bigcup_{i=0}^\infty
F^{-i}(P\cup Z_1)}.
\]

Note, that the sets $A=\overline{Z_3\cup\bigcup_{i=0}^\infty
F^{-i}(P)}$ and $B=\overline{Z_3\cup\bigcup_{i=0}^\infty
F^{-i}(P\cup Z_1)}$ are almost invariant with respect to $F$, in the
sense that
\[
F^{-1}(A)\cup Z_2=A,\quad F^{-1}(B)\cup Z_1\cup Z_2=B,
\]
which is an analog of~\cite{grigorch_z:basilica_sp}*{Theorem~4.1}
for the Basilica group.

\begin{figure}
\centering
\includegraphics[height=170pt]{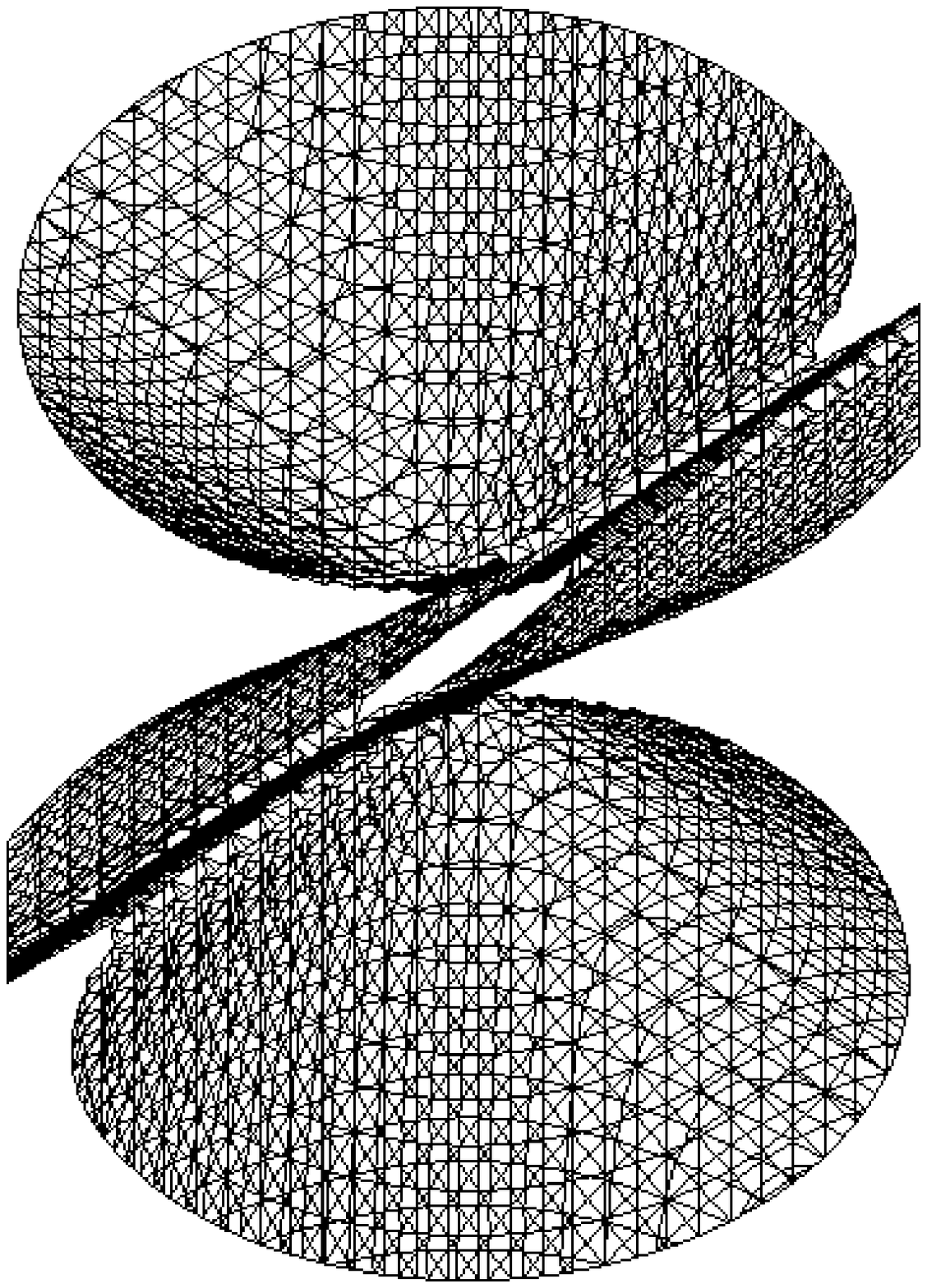}\qquad
\includegraphics[height=170pt]{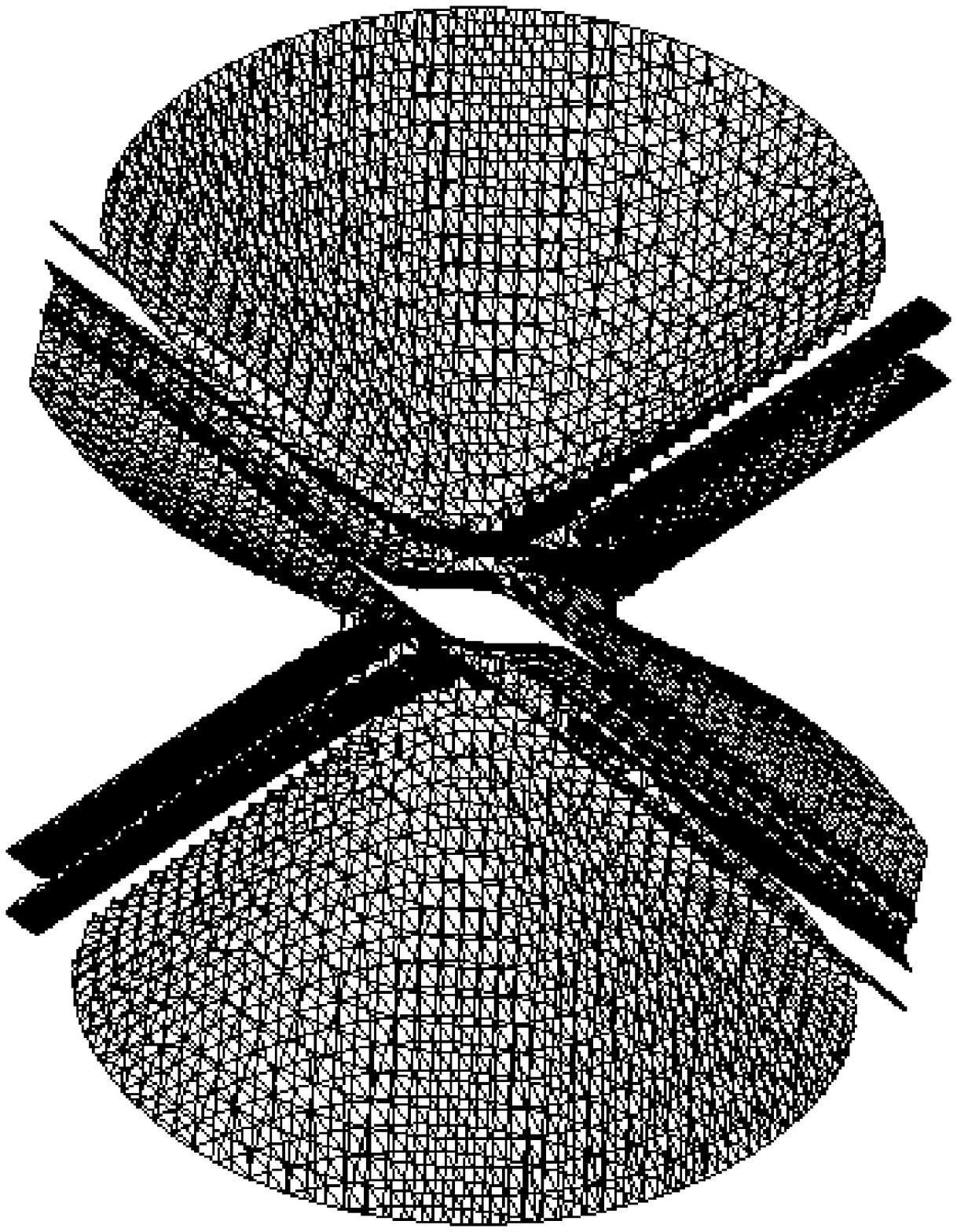}
\caption{Part of the spectrum of
$\widetilde{M}(y,z,\lambda)$\label{fig_preim}}
\end{figure}

The preimages of the plane $P$ under $F^4$ and $F^5$ are shown in
Figure~\ref{fig_preim}.

Note that there are points in the spectrum of $\tilde
M(y,z,\lambda)$ which do not belong to any preimage of the plane
$P$. In particular, the point $(-\frac12,0,\frac12)$ belongs to
$Z_1$ so it is not in the domain of $F$, but
\[
\det\tilde
M_n(-\tfrac12,0,\tfrac12)=\det((a_{n-1}+1)(c_{n-1}+1)-4)=0
\]
since $4$ is an eigenvalue of $(a_{n-1}+1)(c_{n-1}+1)$. However,
this point could be in the closure of the union of all preimages of
$P$.

On the other hand we can formulate a conjecture that the spectrum of
$M=\frac13(a+b+c)$ is the intersection of the line $y=z=1$ with
$A=\overline{Z_3\cup\bigcup_{i=0}^\infty F^{-i}(P)}$, shrunk by a
factor of $\frac13$. This conjecture survives at least up to the
$6$th level.

Note also that the map $F$ is conjugate to a simpler map
\[
G\colon (y,z,\lambda)\to\biggl(\frac zy, \frac\lambda y(-2+y\lambda),
\frac1\lambda(-y+y\lambda^2-\lambda)\biggr)
\]
by the conjugator map
\[
(y,z,\lambda)\to\biggl(\frac1y,\frac1z,y+z-\lambda\biggr).
\]

The histogram for the spectral density of the operator $M_n$ acting
on $9$th level is shown in Figure~\ref{fig_spec}.

Further steps are required to identify the spectrum of the pencil
$\widetilde{M}(y,z,\lambda)$ and of $M$ more precisely. This is related
to the problem of finding invariant subsets of the rational map $F$.
Perhaps the spectrum of $M$ is just the intersection of the
``strange attractor'' of $F$ with the line $y=z=1$, shrunk by factor
of $\frac13$. In any case here we have one more example when the
spectral problem is related to the dynamics of a multidimensional
rational map. There is a hope that the methods developed for this
type of transformations (see, for instance~\cite{sibony:dynamics})
could help to handle this case.


\begin{bibdiv}
\begin{biblist}

\bib{bartholdi_g:spectrum}{article}{%
author={Bartholdi, L.}, author={Grigorchuk, R. I.}, title={On the
spectrum of Hecke type operators related to some fractal
  groups},
journal={trmis},
volume={231},
pages={5--45},
date={2000},
language={Russian},
translation={%
journal={Proc. Steklov Inst. Math.},
volume={2000},
number={4(231)},
pages={1\dash 41}
}
}

\bib{bartholdi_gn:fractal}{article}{%
author={Bartholdi, L.},
author={Grigorchuk, R. I.},
author={Nekrashevych, V.},
title={From fractal groups to fractal sets},
book={%
title={Fractals in Graz 2001},
series={Trends Math.},
pages={25--118},
publisher={Birkh\"auser},
address={Basel},
date={2003}
}
}

\bib{bar_gs:branch}{article}{%
author={Bartholdi, L.}, author={Grigorchuk, R. I.}, author={{\v
S}uni{\'k}, Z.}, title={Branch groups},
book={%
title={Handbook of Algebra, Vol. 3},
pages={989--1112},
publisher={North-Holland},
address={Amsterdam},
date={2003}
}
}

\bib{bartholdi-n:rabbit}{article}{%
author={Bartholdi, L.},
author={Nekrashevych, V.},
title={Thurston equivalence of topological polynomials},
journal={actam1},
volume={197},
date={2006},
pages={1\dash 51}
eprint={arXiv:math.DS/0510082}
}

\bib{bartholdi_v:amenab}{article}{%
author={Bartholdi, L.},
author={Vir{\'a}g, B.},
title={Amenability via random walks},
journal={dukmj},
volume={130},
number={1},
date={2005},
pages={39\dash 56},
eprint={arXiv:math.GR/0305262}
}

\bib{bux_p:iter_monodromy}{article}{%
author={Bux, K.-U.},
author={P{\'e}rez, R.},
title={On the growth of iterated monodromy groups},
book={%
title={Topological and Asymptotic Aspects of Group Theory},
series={Contemp. Math.}, volume={394}, publisher={Amer. Math. Soc.},
address={Providence, RI}, date={2006} }, pages={61\dash 76},
eprint={arXiv:math.GR/0405456}, }


\bib{grigorch:burnside}{article}{%
author={Grigorchuk, R. I.},
title={On {B}urnside's problem on periodic groups},
journal={funap1},
volume={14},
number={1},
pages={53\dash 54},
date={1980},
language={Russian},
translation={%
journal={Functional Anal. Appl.},
volume={14},
number={1},
date={1980},
pages={41\dash 43}
}
}

\bib{grigorch:degrees}{article}{%
author={Grigorchuk, R. I.},
title={Degrees of growth of finitely generated groups and the theory of
  invariant means},
journal={izvan4},
volume={48},
number={5},
pages={939\dash 985},
date={1984},
language={Russian},
translation={%
journal={Math. USSR-Izv.},
volume={25},
number={2},
date={1985},
pages={259\dash 300}
}
}

\bib{grigorch:example}{article}{%
author={Grigorchuk, R. I.},
title={An example of a finitely presented amenable group that does not
  belong to the class EG},
journal={matsb},
volume={189},
number={1},
pages={79\dash100},
date={1998},
language={Russian},
translation={%
journal={Sb. Math.},
volume={189},
number={1-2},
date={1998},
pages={75\dash 95}
}
}

\bib{grigorch:jibranch}{article}{%
author={Grigorchuk, R. I.},
title={Just infinite branch groups},
book={%
title={New Horizons in Pro-$p$ Groups},
volume={184},
series={Progr. Math.},
publisher={Birkh\"auser},
address={Boston, MA},
date={2000}
},
pages={121\dash 179}
}

\bib{gns00:automata}{article}{%
author={Grigorchuk, R. I.},
author={Nekrashevych, V.},
author={Sushchanski{\u\i}, V. I.},
title={Automata, dynamical systems, and groups},
journal={trmis},
volume={231},
pages={134--214},
date={2000},
language={Russian},
translation={%
journal={Proc. Steklov Inst. Math.},
volume={2000},
number={4(231)},
pages={128\dash 203}
}
}

\bib{grigorch_nv:spectra}{article}{%
author={Grigorchuk, R. I.},
author={Nekrashevych, V.},
author={Vorobets, Ya.},
title={Dynamical systems and self-similar groups},
status={preprint},
date={2006}
}

\bib{grigorchuk-s:cr-hanoi}{article}{%
author={Grigorchuk, R. I.},
author={{\v S}uni{\'k}, Z.},
title={Asymptotic aspects of Schreier graphs and Hanoi Towers groups},
journal={crmasp},
volume={342},
number={8},
eprint={arXiv:math.GR/0601592},
date={2006},
pages={545\dash 550}
}

\bib{grigorch_z:asympt_spectrum}{article}{%
author={Grigorchuk, R. I.},
author={{\.Z}uk, A.},
title={On the asymptotic spectrum of random walks on infinite families of
  graphs},
conference={%
title={Random Walks and Discrete Potential Theory},
address={Cortona},
date={1997}
},
book={%
series={Sympos. Math.},
volume={39},
publisher={Cambridge Univ. Press},
address={Cambridge},
date={1999}
},
pages={188\dash 204},
}

\bib{grigorch_z:basilica}{article}{%
author={Grigorchuk, R. I.},
author={{\.Z}uk, A.},
title={On a torsion-free weakly branch group defined by a three state
  automaton},
journal={intjac},
volume={12},
number={1-2},
pages={223\dash 246},
date={2002}
}

\bib{grigorch_z:basilica_sp}{article}{%
author={Grigorchuk, R. I.},
author={{\.Z}uk, A.},
title={Spectral properties of a torsion-free weakly branch group
  defined by a three state automaton},
conference={%
title={Computational and Statistical Group Theory},
address={Las Vegas, NV/Hoboken, NJ},
date={2001}
},
book={%
volume={298},
series={Contemp. Math.},
publisher={Amer. Math. Soc.},
address={Providence, RI},
date={2002}
},
pages={57\dash 82}
}

\bib{kaiman:munchhausen}{article}{%
author={Kaimanovich, V. A.},
title={``M{\"u}nchhausen trick'' and amenability of self-similar
  groups},
journal={intjac},
volume={15},
number={5-6},
pages={907\dash 937},
date={2005}
}

\bib{nekrash:self-similar}{book}{%
author={Nekrashevych, V.},
title={Self-similar groups},
volume={117},
series={Math. Surveys Monogr.},
publisher={ams},
date={2005}
}

\bib{nekrashevych:cantor}{article}{%
author={Nekrashevych, V.}, title={A minimal Cantor set in the space
of $3$-generated groups}, status={to appear in Geom.~Dedicata} }

\bib{sibony:dynamics}{article}{%
author={Sibony, N.},
title={Dynamique des applications rationnelles de $\mathbf{P}^k$},
conference={%
title={Dynamique et g\'eom\'etrie complexes},
address={Lyon},
date={1997}
},
book={%
volume={8},
series={Panor. Synth\`eses},
publisher={Soc. Math. France},
address={Paris},
date={1999},
},
pages={97\dash 185},
}

\end{biblist}
\end{bibdiv}
\end{document}